\documentclass[a4paper,reqno,12pt]{amsart}

\usepackage{graphicx}
\usepackage{mathrsfs}
\usepackage{amsfonts}
\usepackage{amssymb}
\usepackage{amsmath}
\usepackage{amsthm}
\usepackage{amscd}
\usepackage{mathtools}
\usepackage[all,2cell]{xy}
\usepackage{xcolor}
\usepackage[pagebackref,colorlinks]{hyperref}
\usepackage{enumerate}
\usepackage[normalem]{ulem}
\usepackage{mathrsfs}
\usepackage{lineno}
\usepackage{bm}
\usepackage{bbm}
\usepackage{comment}

\numberwithin{equation}{section}

\theoremstyle{plain}
\newtheorem{thm}{Theorem}[section]
\newtheorem{theorem}[thm]{Theorem}
\newtheorem{prop}[thm]{Proposition}
\newtheorem{proposition}[thm]{Proposition}
\newtheorem{cor}[thm]{Corollary}
\newtheorem{lemma}[thm]{Lemma}

\theoremstyle{definition}
\newtheorem{deff}[thm]{Definition}
\newtheorem{definition}[thm]{Definition}
\newtheorem{example}[thm]{Example}
\newtheorem{remark}[thm]{Remark}

\newcommand{\ran}{\bm{r}}
\newcommand{\domr}{\bm{d}}
\newcommand{\G}{\Gamma}

\newcommand{\im}{\operatorname{Im}}
\newcommand{\dom}{\operatorname{Dom}}

\newcommand{\so}{\mathbf{s}}
\newcommand{\ra}{\mathbf{r}}

\newcommand{\gr}{\operatorname{gr}}

\newcommand{\SG}{\operatorname{\mathsf{S}}}
\newcommand{\GP}{\operatorname{\mathsf{G}}}
\newcommand{\inv}{^{-1}}
\def\a{\alpha}

\def\SS{\mathcal{S}}

\def\N{\mathbb N}

\def \Z{\mathbb Z}
\def\-{\text{-}}

\newcommand{\Gc}[1]{\operatorname{\bf Gr-}\boldsymbol #1}

\title[Zhang twists for semigroups and groupoids]{Zhang twists for semigroups and groupoids,\\ with applications to Steinberg algebras}

\author{Roozbeh Hazrat}
\address{R.\ Hazrat: Centre for Research in Mathematics and Data Science, Western Sydney University, AUSTRALIA}
\email{\textcolor[rgb]{0.00,0.00,0.84}{r.hazrat@westernsydney.edu.au}}

\author{Zachary Mesyan}
\address{Z.\ Mesyan: Department of Mathematics, University of Colorado, Colorado Springs, CO, 80918, USA} 
\email{\textcolor[rgb]{0.00,0.00,0.84}{zmesyan@uccs.edu}}

\author{Benjamin Steinberg}
\address{B.\ Steinberg: Department of Mathematics, City College of New York, Convent Avenue at 138th Street, New York, NY 10031, USA} 
\email{\textcolor[rgb]{0.00,0.00,0.84}{bsteinberg@ccny.cuny.edu}}

\subjclass[2020]{20M10, 20M18, 20M30, 22A22, 16S36}

\keywords{Zhang twist, graded semigroup, inverse semigroup, group extension, ample groupoid, graded groupoid, Stone duality, Steinberg algebra, Leavitt path algebra}

\date{\today}

\begin{document}

\begin{abstract}
Zhang twists provide a powerful mechanism for deforming graded algebras, while preserving their graded module categories. We extend this theory by introducing Zhang twists for graded semigroups and graded ample groupoids. We show that the twisting construction preserves a broad class of semigroup structures, and that two graded monoids are twists of each other if and only if their categories of graded sets are pointedly equivalent. We show that Zhang twists are compatible with noncommutative Stone duality: twists on graded-Boolean inverse semigroups are in correspondence with those on the associated graded ample groupoids. Moreover, in the setting of Steinberg algebras we prove that twisting the groupoid gives the same result as twisting the corresponding algebra.

These results provide a unified framework linking graded semigroups, inverse semigroup duality, ample groupoids, and graded algebras, while extending the scope of Zhang twists to settings arising from symbolic dynamics, graph algebras, and Steinberg algebras.  
\end{abstract}

\maketitle

\section{Introduction}

The idea of deforming a mathematical object, while retaining essential features of its structure, has played a fundamental role throughout mathematics. In the setting of associative algebras, this philosophy was developed systematically by Gerstenhaber in his foundational work on deformations of rings and algebras~\cite{GerstenhaberDeformationSeries}. A particularly useful form of deformation for graded algebras arose in the work of Artin, Tate, and Van den Bergh, where the multiplication of a graded algebra is modified using a graded automorphism~\cite{artin}. Zhang subsequently placed this construction in a general framework, by introducing twisting systems for graded algebras~\cite{Zhang}. The resulting Zhang twists have proved to be an important tool in noncommutative algebra and noncommutative geometry, in part because this form of twisting preserves the category of graded modules. For a nice introduction to the relationship between Zhang twists and noncommutative projective geometry see \cite{NamSrivastavaVien2024}. 

The purpose of this paper is to develop a theory of Zhang twists for graded semigroups and graded ample groupoids, and to investigate the interaction of these constructions with noncommutative Stone duality, as well as associated rings, such as Steinberg algebras. Our starting point is the observation that the defining identity for a Zhang twisting system is essentially multiplicative, and therefore makes sense at the level of semigroups. Given a group $\Gamma$ and a $\Gamma$-graded semigroup $S$, we say that a collection $\tau=\{\tau_\alpha\mid \alpha\in\Gamma\}$ of graded bijections of $S$, satisfying an appropriate compatibility condition, is a \emph{Zhang twisting system} on $S$. The multiplication
\[s*t=s\tau_\alpha(t),\]
where $\alpha \in \Gamma$ is the degree of $s$, then produces a new $\Gamma$-graded \emph{twisted} semigroup $S^\tau$. This construction recovers the underlying multiplicative semigroup of the usual Zhang twist of a graded algebra, as defined in~\cite{Zhang}.

A substantial part of the paper is devoted to understanding which semigroup-theoretic properties are preserved under twisting. We show that regularity, being an inverse semigroup, being a graded-Boolean inverse semigroup, and several other fundamental properties are invariant under Zhang twisting (Theorem~\ref{twist-invar}). We also show that a fundamental feature of Zhang twists on algebras, namely the preservation of graded module categories, extends naturally to the semigroup setting (Theorem~\ref{thmprehd6}). Specifically, for graded semigroups $S$ and $T$, 
\[S \cong T^\tau \text{ if and only if } \Gc{S}\cong \Gc{T} \text{ as pointed categories}
\]
(of graded sets). These results allow us to apply the construction in settings where inverse semigroups encode combinatorial, dynamical, or algebraic information arising from discrete structures, such as graphs. 

We further relate Zhang twists on groups to the classical theory of group extensions. A group extension
\[1\longrightarrow N\longrightarrow G\longrightarrow\Gamma\longrightarrow 1\]
gives a $\Gamma$-grading of $G$. We show that, after fixing the induced outer action of $\Gamma$ on $N$, all the group extensions with this outer action can be realized as Zhang twists of one another, by twisting systems acting trivially on $N$ (Theorem~\ref{group-ext}). Thus Zhang twisting provides a semigroup-theoretic interpretation of the familiar cocycle description of group extensions.

We then give an analogous definition of the Zhang twist for graded ample groupoids, and show that twisting preserves the property of being an ample groupoid (Theorem~\ref{gp-twist-thrm}). One of the main motivations for introducing Zhang twists for semigroups and groupoids comes from (graded) noncommutative Stone duality~\cite{Lawson3,RZ2}, which establishes a correspondence between suitable Boolean inverse semigroups and ample groupoids. We show that twisting is compatible with this duality. More precisely, semi-algebraic Zhang twists on graded-Boolean inverse $\land$-semigroups correspond naturally to Zhang twists on their associated ample groupoids, and conversely (Theorem~\ref{duality-thrm}). Moreover, forming the associated groupoid commutes with twisting. This places Zhang deformation within the correspondence between inverse semigroups and groupoids, and shows that the construction is intrinsic to both sides of the duality.

Finally, we examine the effect of twisting an ample groupoid on the associated Steinberg algebra. If $\mathcal{G}$ is a graded ample groupoid, then any  Zhang twist $\tau$ on $\mathcal G$ induces a Zhang twist on the Steinberg algebra $K\mathcal{G}$. We prove that the Steinberg algebra of the twisted groupoid coincides with the corresponding twisting of the original algebra (Theorem~\ref{groupoid-algebra-thrm}):
\[K(\mathcal{G}^\tau) = (K\mathcal{G})^\tau.\]
Consequently, the semigroup and groupoid constructions developed here provide a geometric and combinatorial mechanism for producing Zhang twists on graded algebras. We conclude by considering twisting systems, arising from automorphisms of graph shifts, on graph groupoids and Leavitt path algebras.

The paper is organized as follows. In Section~\ref{semigroupsec} we introduce Zhang twists on graded semigroups, develop their basic properties, study them in various types of semigroups, and relate twists on groups to group extensions. We also introduce twisted semigroup acts (or sets), and establish an analogue of Zhang's equivalence theorem in that context. In Section~\ref{groupoidsec} we define Zhang twists on graded ample groupoids, and prove their compatibility with graded noncommutative Stone duality. Section~\ref{algebrasec} is devoted to Zhang twists on various rings, particularly corner skew Laurent polynomial rings and Steinberg algebras, and applies the theory to graph groupoids, dynamical systems arising from shifts, and Leavitt path algebras.

\section{Semigroups}\label{semigroupsec}

\subsection{Graded semigroups} 

We begin by recalling semigroup terminology that will be used extensively.

A \emph{semigroup} is a nonempty set equipped with an associative binary operation. A \emph{monoid} is a semigroup with an identity element $1$.  We assume throughout that each semigroup has a zero element $0$, unless specified otherwise. A semigroup $S$ is \emph{regular} if for each $s\in S$ there exists a $t\in S$ such that $sts=s$. One can show that $S$ is regular if and only if for every $s\in  S$ there is an \emph{inner inverse} $t\in S$ satisfying $sts=s$ and $tst=t$. If every $s \in S$ has a unique inner inverse, denoted $s^{-1}$, then $S$ is called an \emph{inverse semigroup}. We refer the reader to~\cite{Howie} for the theory of semigroups, and to~\cite{Lawson1} for that of inverse semigroups. 
 
Let $S$ be a semigroup and $\Gamma$ a group. Then $S$ is $\Gamma$-\emph{graded} if there is a function $\phi : S \setminus\{0\} \to \Gamma$ such that $\phi(st) = \phi(s)\phi(t)$, whenever $st \not= 0$. We set $S_\alpha:=\phi^{-1}(\alpha) \cup \{0\}$ for each $\alpha \in \Gamma$. Equivalently, $S$ is $\Gamma$-\emph{graded} if $S = \bigcup_{\alpha \in \Gamma} S_\alpha$ for some $S_{\alpha} \subseteq S$, such that $S_{\alpha} S_{\beta} \subseteq S_{\alpha\beta}$ for all $\alpha,\beta \in \Gamma$, and $S_{\alpha} \cap S_{\beta} = \{0\}$ for all distinct $\alpha,\beta \in \Gamma$. If $S$ is $\Gamma$-graded, then we say that it is \emph{strongly $\Gamma$-graded} if $S_{\alpha} S_{\beta} = S_{\alpha\beta}$ for all $\alpha,\beta \in \Gamma$.

Let $S$ be a $\Gamma$-graded semigroup. For each $\alpha \in \Gamma$ and $s \in S_{\alpha} \setminus \{0\}$, we say that the \emph{degree} of $s$ is $\alpha$, and write $|s|:=\alpha$. We refer to $S_{\alpha}$ as the \emph{homogeneous component of $S$ of degree $\alpha$}, for each $\alpha \in \Gamma$. Letting $E(S)$ denote the set of idempotents in $S$ and $\varepsilon$ the identity element of $\Gamma$, we note that $|s|=\varepsilon$ for all $s \in E(S) \setminus \{0\}$. Moreover, if $S$ is an inverse semigroup, then $|s^{-1}| = |s|^{-1}$ for all $s \in S \setminus \{0\}$. The set $\big\{\alpha \in \Gamma \mid S_\alpha \not = \{0\} \big\}$ is called the \emph{support of $S$}. 

Given $\Gamma$-graded semigroups $S$ and $T$, a function $\phi:S \rightarrow T$ is  said to be \emph{graded} if $\phi(S_\alpha)\subseteq T_\alpha$ for every $\alpha \in \Gamma$. See~\cite{RZ1} for more details about graded semigroups.

\subsection{Twisted semigroups}

We are now ready to define our main object of study, which is based on the concept of twisting systems for algebras given in~\cite[Definition 2.1]{Zhang}. (See Example~\ref{ring-twist} below for more details.)

\begin{definition} \label{twist-def}
Let $\Gamma$ be a group and $S$ a $\Gamma$-graded semigroup. A \emph{Zhang twisting system}, or a \emph{Zhang twist, $\tau$ on $S$}, is a collection $\{\tau_{\alpha} \mid \alpha \in \Gamma\}$, where $\tau_{\alpha} : S \to S$ is a graded bijection for each $\alpha \in \Gamma$, and
\begin{equation} \label{twist-eq}
\tau_{\alpha}(s\tau_{\beta}(t)) = \tau_{\alpha}(s)\tau_{\alpha \beta}(t)
\end{equation}
for all $\alpha, \beta \in \Gamma$, $s \in S_{\beta}$, and $t \in S$. A Zhang twist $\tau = \{\tau_{\alpha} \mid \alpha \in \Gamma\}$ on $S$ is \emph{algebraic} if $\tau_{\alpha}$ is a semigroup automorphism and $\tau_{\alpha}\tau_{\beta} = \tau_{\alpha\beta}$, for all $\alpha, \beta \in \Gamma$.
\end{definition}

A Zhang twist on a graded semigroup allows one to define an alternative multiplication operation.

\begin{lemma} \label{twist-mult}
Let $\, \Gamma$ be a group, $\, (S, \cdot)$ a semigroup $\, \Gamma$-graded by $\phi : S \setminus\{0\} \to \Gamma$, and $\tau = \{\tau_{\alpha} \mid \alpha \in \Gamma\}$ a Zhang twist on $S$. Define a binary operation $*$ on $S$ via 
\[s * t := s \tau_{\alpha}(t),\] 
for all $\alpha \in \Gamma$, $s \in S_{\alpha}$, and $t \in S$. Then $\, (S,*)$ is a semigroup, $\, \Gamma$-graded by $\phi$.
\end{lemma}

\begin{proof}
Let $\alpha, \beta, \gamma \in \Gamma$, $r \in S_{\alpha}$, $s \in S_{\beta}$, and $t \in S_{\gamma}$. Since $\tau_{\alpha}(s) \in S_{\beta}$,
\[r*(s*t) = r\tau_{\alpha}(s\tau_{\beta}(t)) = r\tau_{\alpha}(s)\tau_{\alpha\beta}(t) = (r*s)*t.\]
Thus $*$ is associative, and hence $(S,*)$ is a semigroup. If $s*t \neq 0$, then $s \tau_{\beta}(t) \neq 0$, and so
\[\phi(s*t) = \phi(s \tau_{\beta}(t)) = \phi(s)\phi(\tau_{\beta}(t)) = \beta\gamma = \phi(s)\phi(t),\]
since $\tau_{\beta}$ is graded. Thus $\phi$ defines a grading on $(S,*)$.
\end{proof}

\begin{definition}
Let $\Gamma$ be a group, $S$ a $\Gamma$-graded semigroup, $\tau$ a Zhang twist on $S$, and $*$ as in Lemma~\ref{twist-mult}. We denote the semigroup $(S,*)$ by $S^{\tau}$, and refer to it as the \emph{(Zhang) twisting of $S$ by $\tau$}, or, more simply, as a \emph{twisted} semigroup.
\end{definition}

\begin{remark}
The above definition of ``twisting" is clearly not left-right symmetric. Analogously to~\cite[Definition 4.1]{Zhang}, one can define a \emph{left Zhang twist} $\nu$ on a $\Gamma$-graded semigroup $S$ to be a collection $\{\nu_{\alpha} \mid \alpha \in \Gamma\}$, where $\nu_{\alpha} : S \to S$ is a graded bijection for each $\alpha \in \Gamma$, and
\[\nu_{\alpha}(\nu_{\beta}(s)t) = \nu_{\beta\alpha}(s)\nu_{\alpha}(t)\]
for all $\alpha, \beta \in \Gamma$, $s \in S$, and $t \in S_{\beta}$. Then one can show that, in this situation, the operation $s * t := \nu_{\alpha}(s)t$ ($\alpha \in \Gamma$, $s \in S$, $t \in S_{\alpha}$), turns $(S,*)$ into a semigroup $\prescript{\nu}{}S$, with the same grading.

It is easy to see that a \emph{left twisting} of $S$, as above, corresponds to a twisting of $S^{\mathrm{op}}$, the opposite semigroup of $S$, graded by $\Gamma^{\mathrm{op}}$, the opposite group of $\Gamma$, and vice versa. So a separate theory of left twistings is unnecessary. Moreover, using the same argument as in~\cite[Theorem 4.3]{Zhang}, one can show that given a Zhang twist $\tau = \{\tau_{\alpha} \mid \alpha \in \Gamma\}$ on $S$, letting $\nu_{\alpha}(s) = \tau_{\beta\alpha}\tau_{\beta}^{-1}(s)$ for all $\alpha, \beta \in \Gamma$ and $s \in S_{\beta}$, defines a left twisting $\nu$ of $S$, such that $\prescript{\nu}{}S \cong S^{\tau}$ as graded semigroups.
\end{remark}

\subsection{Examples} 

Before developing a theory of twisted semigroups, let us discuss ways to construct them. The following two observations are analogues of ones from~\cite{Zhang}.

\begin{remark}\label{alg-eg}
Let $\Gamma$ be a group, $S$ a $\Gamma$-graded semigroup, and $\mathrm{Aut}_{\mathrm{gr}}(S)$ the group of graded semigroup automorphisms of $S$. Suppose that $\tau : \Gamma \rightarrow \mathrm{Aut}_{\mathrm{gr}}(S)$ is a homomorphism, and let $\tau_{\alpha} := \tau(\alpha)$ for each $\alpha \in \Gamma$. Then for all $\alpha, \beta \in \Gamma$, $s \in S_{\beta}$, and $t \in S$ we have 
\[\tau_{\alpha}(s\tau_{\beta}(t)) = \tau_{\alpha}(s)\tau_{\alpha}(\tau_{\beta}(t)) = \tau_{\alpha}(s)\tau_{\alpha \beta}(t),\]
since $\tau(\alpha)\tau(\beta) = \tau(\alpha\beta)$. Therefore $\{\tau_{\alpha} \mid \alpha \in \Gamma\}$ is an algebraic Zhang twist on $S$. 
\end{remark}

\begin{remark}\label{aut-eg}
Let $S$ be a $\Z$-graded semigroup, where $\Z$ denotes the group of the integers, and let $\psi : S \rightarrow S$ be a graded automorphism. Then $\tau(n) := \psi^n$ gives a homomorphism $\tau : \Z \rightarrow \mathrm{Aut}_{\mathrm{gr}}(S)$, and so $\{\psi^n \mid n \in \Z\}$ is an algebraic Zhang twist on $S$, by Remark~\ref{alg-eg}.
\end{remark}

Next we give a couple of examples of $\Z$-graded semigroups, along with graded automorphisms thereof. Additional $\Z$-graded semigroups can be found in~\cite[Examples 2.3, 2.4, 2.8, 2.11, 2.18]{RZ1}.

\begin{example}
Let $X$ be a set, and let $F=\langle X \rangle$ denote either the free semigroup or the free commutative semigroup (with zero) generated by $X$. Letting $F_n$ denote the subset of $F$ consisting of words of length $n$ in the generators, together with $0$, we obtain a $\Z$-grading $F=\bigcup_{i\geq 0} F_n$, with support $\N$ (the set of the natural numbers). Any bijection $f : X \to X$ can be extended to a graded automorphism $\psi_f : F \rightarrow F$, via $\psi_f(0) := 0$, and $\psi_f(x_{1} \cdots x_{n}) := f(x_{1}) \cdots f(x_{n})$ for all $x_{1}, \dots, x_{n} \in X$. Then $\psi_f$ can be used to construct a Zhang twist on $F$, as in Remark~\ref{aut-eg}.

Now consider the particular case where $F$ is the free commutative semigroup generated by $\{x,y\}$. Let $f : \{x,y\} \rightarrow \{x,y\}$ be defined by $f(x) := y$ and $f(y) := x$, extend $f$ to an automorphism $\psi_f : F \rightarrow F$, as above, and let $\tau := \{\psi_f^n \mid n \in \Z\}$ be the resulting Zhang twist on $F$. Then in $F^{\tau}$ we have $x * y = x \psi_f(y) = x^2$ and $y * x = y \psi_f(x) = y^2$. So $F^{\tau}$ is not commutative, and therefore not isomorphic to $F$. It is also not hard to see that $F^{\tau}$ is not isomorphic to the free semigroup on two generators.
\end{example}

\begin{example}\label{graphinsemigrp}
Recall that a \emph{directed graph} $E=(E^0,E^1,\ra,\so)$ consists of two disjoint sets, $E^0$ and $E^1$ (containing \emph{vertices} and \emph{edges}, respectively), together with functions $\so,\ra : E^1 \to E^0$, called \emph{source} and \emph{range}, respectively. A \emph{path} $x$ in $E$ is a finite sequence of edges $x=e_1\cdots e_n$ such that $\ra(e_i)=\so(e_{i+1})$ for all $i \in \{1,\dots,n-1\}$. In this case, $\so(x):=\so(e_1)$ is the \emph{source} of $x$, $\ra(x):=\ra(e_n)$ is the \emph{range} of $x$, and $|x|:=n$ is the \emph{length} of $x$. We view vertices as paths of length $0$, and extend $\so$ and $\ra$ to $E^0$ by setting $\so(v) : = v$ and $\ra(v) := v$ for all $v \in E^0$. We refer the reader to \cite{Abrams-Monograph} for background on directed graphs and associated algebras.

Given a directed graph $E=(E^0,E^1,\ra,\so)$, the \emph{grqaph inverse semigroup $\SS(E)$ of $E$} is the semigroup (with zero) generated by $E^0 \cup E^1 \cup \{e^{-1} \mid e\in E^1\}$, subject to the following relations for all $v,w\in E^0$ and $e,f\in E^1$ (with $\delta$ denoting the Kronecker delta):
\begin{enumerate}
\item[(V)] $vw = \delta_{v,w}v$;
\item[(E1)] $\so(e)e=e\ra(e)=e$;
\item[(E2)] $\ra(e)e^{-1}=e^{-1}\so(e)=e^{-1}$;
\item[(CK1)] $e^{-1}f=\delta _{e,f}\ra(e)$.
\end{enumerate}
We define $v^{-1}=v$ for each $v \in E^0$, and for any path $y=e_1\cdots e_n$ ($e_1,\dots, e_n \in E^1$) we let $y^{-1} = e_n^{-1} \cdots e_1^{-1}$. With this notation, it is easy to see that every nonzero element of $\SS(E)$ can be written as $xy^{-1}$, for some paths $x$ and $y$ satisfying $\ra(x)=\ra(y)$, and $\SS(E)$ is indeed an inverse semigroup, with $(xy^{-1})^{-1} = yx^{-1}$. Moreover $\SS(E)$ is naturally $\Z$-graded via $\phi : \SS(E) \setminus \{0\} \to \Z$, where $\phi(xy^{-1}) := |x|-|y|$, for all paths $x$ and $y$ in $E$.

A \emph{graph automorphism} of $E$ is a pair of bijections $\phi_0 : E^0 \to E^0$ and $\phi_1 : E^1 \to E^1$, such that $\phi_0(\so(e)) = \so(\phi_1(e))$ and $\phi_0(\ra(e)) = \ra(\phi_1(e))$, for every $e \in E^1$. It is not hard to show (see, e.g.,~\cite[Theorem 20]{MM}) that any automorphism of $E$ can be extended to a, necessarily graded, automorphism of $\SS(E)$, and hence also to a Zhang twist, as in Remark~\ref{aut-eg}.
\end{example}

The following is an example of a Zhang twist that does not result from an automorphism of a $\Z$-graded semigroup.

\begin{example}
Let $X$ be a nonempty set. Recall that a \emph{partial symmetry} of $X$ is a function $\phi: A\rightarrow B$, where $A,B \subseteq X$. Here we let $\dom(\phi):=A$ and $\im(\phi): = B$. The set $\mathcal{I}(X)$ of all partial symmetries of $X$ is an inverse semigroup with respect to composition of relations, called the \emph{symmetric inverse monoid}. Specifically, for all $\phi,\psi \in \mathcal{I}(X)$ we take $\phi \psi$ to be the composite of $\phi$ and $\psi$ as functions, restricted to the domain $\psi^{-1}(\im(\psi)\cap \dom(\phi))$. The empty function plays the role of the zero element in $\mathcal{I}(X)$.  

Next, let $\Gamma$ be a group, and let $\psi : X \rightarrow \Gamma$ a function. For each $\alpha \in \Gamma$ and nonempty $Y \subseteq X$, let $Y_{\alpha}:= Y\cap \psi^{-1}(\alpha)$ and
\begin{equation*}\label{mothercomp}
\mathcal{I}(X)_{\alpha}:= \big \{\phi \in \mathcal{I}(X) \, \big \vert  \, \phi(\dom(\phi)_{\beta}) \subseteq X_{\alpha\beta} \text{ for all } \beta \in \Gamma \big \}. 
\end{equation*}
Then
\begin{equation*}\label{mothersemigroup}
\mathcal{I}^{\gr}(X):= \bigcup _{\alpha \in \Gamma} \mathcal{I}(X)_{\alpha} 
\end{equation*}
is a $\Gamma$-graded inverse semigroup, according to~\cite[Proposition 5.3]{RZ1}.

Suppose that $\rho : \Gamma \rightarrow \mathrm{Sym}(X)$ is a homomorphism such that $\rho(\alpha)(X_{\beta}) = X_{\beta}$ for all $\alpha, \beta \in \Gamma$, where $\mathrm{Sym}(X)$ denotes the group of all permutations of $X$. Then letting $\tau(\alpha)(\phi) := \rho(\alpha)\phi \rho(\alpha)^{-1}$, for each $\alpha \in \Gamma$ and $\phi \in \mathcal{I}^{\gr}(X)$, defines a homomorphism $\tau: \Gamma \rightarrow \mathrm{Aut}_{\mathrm{gr}}(\mathcal{I}^{\gr}(X))$, which leads to an algebraic Zhang twist on $\mathcal{I}^{\gr}(X)$, as in Remark~\ref{alg-eg}. 

To construct a concrete example of this sort, take $\Gamma$ to be any group, let $X := \bigsqcup_{\alpha \in \Gamma} X_{\alpha}$ be the disjoint union, where each $X_{\alpha}$ is a copy of $\Gamma$, and define $\rho(\alpha)(\gamma) := \alpha\gamma \in X_{\beta}$, for all $\alpha, \beta \in \Gamma$ and $\gamma \in X_{\beta}$.
\end{example}

\begin{example} \label{ring-twist}
Recall that given a group $\Gamma$, a ring $A$ is $\Gamma$-\emph{graded} if $A = \bigoplus_{\alpha \in \Gamma} A_{\alpha}$, where the $A_{\alpha}$ are additive subgroups of $A$, and $A_{\alpha} A_{\beta} \subseteq A_{\alpha \beta}$ for all $\alpha, \beta \in \Gamma$ (with $A_{\alpha} A_{\beta}$ denoting the set of all sums of elements of the form $ab$, with $a \in A_{\alpha}$ and $b \in A_{\beta}$). For any $\Gamma$-graded ring $A$, letting $S:= \bigcup_{\alpha \in \Gamma} A_{\alpha}$ gives a $\Gamma$-graded semigroup, under multiplication. 

Now let $K$ be a commutative unital ring, $\Gamma$ a group, and $A$ a $\Gamma$-graded $K$-algebra. According to~\cite[Definition 2.1]{Zhang}, a \emph{(Zhang) twisting system $\tau$ on $A$} is a collection $\{\tau_{\alpha} \mid \alpha \in \Gamma\}$, where $\tau_{\alpha} : A \to A$ is a $K$-linear graded bijection for each $\alpha \in \Gamma$, and
\[\tau_{\alpha}(a\tau_{\beta}(b)) = \tau_{\alpha}(a)\tau_{\alpha \beta}(b)\]
for all $\alpha, \beta, \gamma \in \Gamma$, $a \in A_{\beta}$, and $b \in A_{\gamma}$. Clearly, in this situation, restricting each $\tau_{\alpha}$ to $S := \bigcup_{\alpha \in \Gamma} A_{\alpha}$ leads to a Zhang twist, in the sense of Definition~\ref{twist-def}. 
\end{example}

\begin{remark}
Let $\Gamma$ be a group, $K$ a commutative ring, $A = \bigoplus_{\alpha \in \Gamma} A_{\alpha}$ a $\Gamma$-graded $K$-algebra, and $S := \bigcup_{\alpha \in \Gamma} A_{\alpha}$ (as in Example~\ref{ring-twist}). Suppose that $\tau = \{\tau_{\alpha} \mid \alpha \in \Gamma\}$ is a Zhang twist on $S$, such that each $\tau_{\alpha} : S \rightarrow S$ extends to a $K$-linear function $\tau_{\alpha} : A \to A$. Then it is easy to see that extending each $\tau_{\alpha}$ in that manner gives a twisting system on $A$, in the sense of~\cite[Definition 2.1]{Zhang}, and every twisting system on $A$ arises this way.
\end{remark}

We give further illustrations in Sections~\ref{cornerskewsemigroup} and~\ref{gp-sect}, including a non-algebraic Zhang twist on a group (Example~\ref{gp-eg}).

\subsection{Twisting basics}

In this subsection we record some basic observations about Zhang twists on semigroups, that will be used extensively. These are mostly analogues of results from~\cite{Zhang}. 

\begin{lemma} \label{id-lem}
Let $\, \Gamma$ be a group, $S$ a $\, \Gamma$-graded semigroup, and $\tau = \{\tau_{\alpha} \mid \alpha \in \Gamma\}$ a Zhang twist on $S$. Then the following hold for all $\alpha, \beta \in \Gamma$, $s \in S_{\beta}$, and $t \in S$:
\begin{enumerate}
\item[$(1)$] $\tau_{\alpha}(st) = \tau_{\alpha}(s)\tau_{\alpha \beta}\tau_{\beta}^{-1}(t)$;
\item[$(2)$] $\tau_{\alpha}^{-1}(s\tau_{\alpha\beta}(t)) = \tau_{\alpha}^{-1}(s)\tau_{\beta}(t)$;
\item[$(3)$] $\tau_{\alpha}^{-1}(st) = \tau_{\alpha}^{-1}(s)\tau_{\beta}\tau_{\alpha \beta}^{-1}(t)$.
\end{enumerate}
\end{lemma}

\begin{proof}
The equation in (1) is the result of substituting $\tau_{\beta}^{-1}(t)$ for $t$ in equation (\ref{twist-eq}) in Definition~\ref{twist-def}. The equation in (2) is the result of substituting $\tau_{\alpha}^{-1}(s)$ for $s$ in equation (\ref{twist-eq}), and then applying $\tau_{\alpha}^{-1}$ to both sides. Finally, the equation in (3) is the result of substituting $\tau_{\alpha\beta}^{-1}(t)$ for $t$ in the equation in (2).
\end{proof}

\begin{lemma} \label{twist-id-elt}
Let $\, \Gamma$ be a group, $S$ a $\, \Gamma$-graded semigroup, and $\tau = \{\tau_{\alpha} \mid \alpha \in \Gamma\}$ a Zhang twist on $S$.
\begin{enumerate}
\item[$(1)$] For all $\alpha \in \Gamma$ we have $\tau_{\alpha}(0) = 0$, and hence $\, 0$ is the zero element of $S^{\tau}$.
\item[$(2)$] If $S$ is a monoid, then $\tau_{\varepsilon}(1) = 1$ if and only if $\tau_{\varepsilon}(s) = s$ for all $s \in S$.
\item[$(3)$] If $S$ is a monoid, then $\, 1 = \tau_{\alpha}^{-1}(1)\tau_{\varepsilon}(1)$ for all $\alpha \in \Gamma$. 
\item[$(4)$] If $S$ is an inverse semigroup, then $\tau_{\varepsilon}(u) = u$ for all $u \in E(S)$  if and only if $\tau_{\varepsilon}(s) = s$ for all $s \in S$.
\end{enumerate}
\end{lemma}

\begin{proof}
(1) The first claim follows immediately from the hypothesis that $\tau_{\alpha}$ is a graded bijection of $S$, for each $\alpha \in \Gamma$, and the fact that $0$ is the unique element of $S$ common to every component $S_{\beta}$. The second claim is an easy consequence of the first.

\smallskip

(2) Suppose that $S$ is a monoid and $\tau_{\varepsilon}(1) = 1$, and let $s \in S$. Then
\[\tau_{\varepsilon}(\tau_{\varepsilon}(s)) = \tau_{\varepsilon}(1\cdot \tau_{\varepsilon}(s)) = \tau_{\varepsilon}(1)\tau_{\varepsilon}(s) = \tau_{\varepsilon}(s).\]
Since $\tau_{\varepsilon}$ is a bijection, this implies that $\tau_{\varepsilon}(s) = s$.

The converse is tautological.

\smallskip

(3) Suppose that $S$ is a monoid, and let $\alpha \in \Gamma$. Then, using (3) in Lemma~\ref{id-lem}, 
\[1 = \tau_{\alpha}^{-1}\tau_{\alpha}(1) = \tau_{\alpha}^{-1}(1 \cdot \tau_{\alpha}(1)) = \tau_{\alpha}^{-1}(1)\tau_{\varepsilon}(\tau_{\alpha\varepsilon}^{-1}\tau_{\alpha}(1)) = \tau_{\alpha}^{-1}(1)\tau_{\varepsilon}(1).\]

\smallskip

(4) Suppose that $S$ is an inverse semigroup, and $\tau_{\varepsilon}(u) = u$ for all $u \in E(S)$. Let $s \in S$, and set $u := \tau_{\varepsilon}(s)\tau_{\varepsilon}(s)^{-1} \in E(S) \subseteq S_{\varepsilon}$. Then
\[\tau_{\varepsilon}(\tau_{\varepsilon}(s)) = \tau_{\varepsilon}(u\tau_{\varepsilon}(s)) =  \tau_{\varepsilon}(u)\tau_{\varepsilon}(s) = u\tau_{\varepsilon}(s) = \tau_{\varepsilon}(s),\]
and so $\tau_{\varepsilon}(s) = s$, given that $\tau_{\varepsilon}$ is a bijection.

The converse is tautological.
\end{proof}

\begin{proposition} \label{id-fix}
Let $\, \Gamma$ be a group, $S$ a $\, \Gamma$-graded semigroup, and $\tau$ a Zhang twist on $S$. Then there is a Zhang twist $\sigma$ on $S$ such that $\sigma_{\varepsilon}(s) = s$ for all $s \in S$, and $S^{\tau} \cong S^{\sigma}$, via a graded isomorphism.
\end{proposition}

\begin{proof}
Write $\tau = \{\tau_{\alpha} \mid \alpha \in \Gamma\}$, choose $\gamma \in \Gamma$, and define $\sigma_{\alpha} := \tau_{\gamma\alpha}\tau_{\gamma}^{-1}$ for each $\alpha \in \Gamma$. Then $\sigma_{\varepsilon} (s) = \tau_{\gamma}\tau_{\gamma}^{-1}(s) = s$ for all $s\in S$, and clearly $\sigma_{\alpha} : S \to S$ is a graded bijection for each $\alpha \in \Gamma$. Now, let $\alpha, \beta \in \Gamma$, $s \in S_{\beta}$, and $t \in S$. Then, using (2) in Lemma~\ref{id-lem}, 
\begin{multline*}
\sigma_{\alpha}(s\sigma_{\beta}(t)) = \tau_{\gamma\alpha}\tau_{\gamma}^{-1}(s\tau_{\gamma\beta}\tau_{\gamma}^{-1}(t)) = \tau_{\gamma\alpha}(\tau_{\gamma}^{-1}(s)\tau_{\beta}\tau_{\gamma}^{-1}(t))\\
= \tau_{\gamma\alpha}\tau_{\gamma}^{-1}(s)\tau_{\gamma\alpha\beta}\tau_{\gamma}^{-1}(t) = \sigma_{\alpha}(s)\sigma_{\alpha\beta}(t).
\end{multline*}
Hence $\sigma = \{\sigma_{\alpha} \mid \alpha \in \Gamma\}$ is a Zhang twist on $S$.

Finally, let $\star$ denote the multiplication in $S^{\sigma}$. Then for all $\alpha \in \Gamma$, $s \in S_{\alpha}$, and $t \in S$,
\[\tau_{\gamma}(s*t) = \tau_{\gamma}(s \tau_{\alpha}(t)) = \tau_{\gamma}(s) \tau_{\gamma\alpha}(t) = \tau_{\gamma}(s) \tau_{\gamma\alpha}\tau_{\gamma}^{-1}\tau_{\gamma}(t) = \tau_{\gamma}(s) \sigma_{\alpha}(\tau_{\gamma}(t)) = \tau_{\gamma}(s) \star \tau_{\gamma}(t).\]
It follows that $\tau_{\gamma} : S^{\tau} \to S^{\sigma}$ is an isomorphism, and it is graded, by Lemma~\ref{twist-mult}.
\end{proof}

\begin{proposition} \label{inverse-twist}
Let $\, \Gamma$ be a group, $S$ a $\, \Gamma$-graded semigroup, and $\tau = \{\tau_{\alpha} \mid \alpha \in \Gamma\}$ a Zhang twist on $S$. 
\begin{enumerate}
\item[$(1)$] The set $\tau^{-1}: = \{\tau_{\alpha}^{-1} \mid \alpha \in \Gamma\}$ is a Zhang twist on $S^{\tau}$.
\item[$(2)$] For any Zhang twist $\sigma = \{\sigma_{\alpha} \mid \alpha \in \Gamma\}$ on $S^{\tau}$, the set $\tau\sigma := \{\tau_{\alpha}\sigma_{\alpha} \mid \alpha \in \Gamma\}$ is a Zhang twist on $S$, and $S^{\tau\sigma} = (S^{\tau})^{\sigma}$, as semigroups.
\end{enumerate}
\end{proposition}

\begin{proof}
(1) Clearly $\tau_{\alpha}^{-1} : S \to S$ is a graded bijection for each $\alpha \in \Gamma$. Also, for all $\alpha, \beta \in \Gamma$, $s \in S_{\beta}$, and $t \in S$, using (3) in Lemma~\ref{id-lem} gives
\[\tau_{\alpha}^{-1}(s * \tau_{\beta}^{-1}(t)) = \tau_{\alpha}^{-1}(st) = \tau_{\alpha}^{-1}(s)\tau_{\beta}\tau_{\alpha \beta}^{-1}(t)  = \tau_{\alpha}^{-1}(s) * \tau_{\alpha \beta}^{-1}(t),\]
where $*$ is the multiplication in $S^{\tau}$. Hence $\tau^{-1}$ is a Zhang twist on $S^{\tau}$.

\smallskip

(2) As above, $\tau_{\alpha}\sigma_{\alpha} : S \to S$ is certainly a graded bijection for each $\alpha \in \Gamma$. Now, for all $\alpha, \beta \in \Gamma$, $s \in S_{\beta}$, and $t \in S$, we have
\begin{multline*}
\tau_{\alpha}\sigma_{\alpha} (s\tau_{\beta}\sigma_{\beta} (t)) = \tau_{\alpha}\sigma_{\alpha} (s*\sigma_{\beta} (t)) = \tau_{\alpha} (\sigma_{\alpha}(s)*\sigma_{\alpha\beta} (t))\\
= \tau_{\alpha} (\sigma_{\alpha}(s)\tau_{\beta}\sigma_{\alpha\beta} (t)) = \tau_{\alpha} \sigma_{\alpha}(s)\tau_{\alpha\beta}\sigma_{\alpha\beta} (t),
\end{multline*}
where $*$ is the multiplication in $S^{\tau},$ and so $\tau\sigma$ is a Zhang twist on $S$. Next, let $\star$, respectively $\bullet$, denote the multiplication in $S^{\tau\sigma}$, respectively $(S^{\tau})^{\sigma}$. Then for all $\alpha \in \Gamma$, $s \in S_{\alpha}$, and $t \in S$,
\[s \star t = s\tau_{\alpha}\sigma_{\alpha}(t) = s*\sigma_{\alpha}(t) = s \bullet t.\]
It follows that $(S^{\tau\sigma}, \star) = ((S^{\tau})^{\sigma},\bullet)$.
\end{proof}

\subsection{Twisted inverse semigroups} 

Twisting preserves much of the structure of an inverse semigroup, as we show next.

\begin{lemma} \label{idempt-lem}
Let $\, \Gamma$ be a group, $S$ a $\, \Gamma$-graded inverse semigroup, and $\tau = \{\tau_{\alpha} \mid \alpha \in \Gamma\}$ a Zhang twist on $S$, such that $\tau_{\varepsilon}(u)=u$ for all $u \in E(S)$. Then $E(S) = \tau_{\alpha}(E(S)) = E(S^{\tau})$ for all $\alpha \in \Gamma$.
\end{lemma}

\begin{proof}
Let $u \in S$ and $\alpha \in \Gamma$. If $u \in E(S)$, then $u \in S_{\varepsilon}$, and so, using the hypothesis on $\tau_{\varepsilon}$, gives
\[u*u = u \tau_{\varepsilon}(u) = u^2 = u\]
and $u \in E(S^{\tau})$. Additionally, if $u \in E(S)$, then
\[u = \tau_{\alpha}^{-1}\tau_{\alpha}(u^2) = \tau_{\alpha}^{-1}\tau_{\alpha}(u\tau_{\varepsilon}(u)) = \tau_{\alpha}^{-1}(\tau_{\alpha}(u)\tau_{\alpha}(u)),\]
which implies that $\tau_{\alpha}(u) = \tau_{\alpha}(u)^2$. Thus $\tau_{\alpha}(E(S)) \subseteq E(S) \subseteq E(S^{\tau})$.

Next, suppose that $\tau_{\alpha}(u) \in E(S)$. Then $\tau_{\alpha}(u) \in S_{\varepsilon}$, and so
\[\tau_{\alpha}(u) = \tau_{\alpha}(u)\tau_{\alpha}(u) = \tau_{\alpha}(u\tau_{\varepsilon}(u)) = \tau_{\alpha}(u^2).\]
Since $\tau_{\alpha}$ is a bijection, this means that $u = u^2$, and so $u \in E(S)$. Since $\tau_{\alpha}(E(S)) \subseteq E(S)$, using once again the fact that $\tau_{\alpha}$ is a bijection on $S$, we conclude that $E(S) = \tau_{\alpha}(E(S))$.

Finally, suppose that $u \in E(S^{\tau})$. Then $u \in S_{\varepsilon}^{\tau} = S_{\varepsilon}$, and so
\[u = u*u = u \tau_{\varepsilon}(u) = u^2.\]
Thus, $u \in E(S)$, and therefore $E(S) = E(S^{\tau})$.
\end{proof}

\begin{proposition} \label{twist-reg}
Let $\, \Gamma$ be a group, $S$ a $\, \Gamma$-graded semigroup, and $\tau = \{\tau_{\alpha} \mid \alpha \in \Gamma\}$ a Zhang twist on $S$. 
\begin{enumerate}
\item[$(1)$] $S$ is regular if and only if $S^{\tau}$ is.
\item[$(2)$] $S$ is an inverse semigroup if and only if $S^{\tau}$ is.
\end{enumerate}
\end{proposition}

\begin{proof}
By Proposition~\ref{inverse-twist}, it suffices to prove only the forward implication of each statement. Moreover, each of the properties in question is preserved by isomorphism, and so, by Proposition~\ref{id-fix}, we may assume that $\tau_{\varepsilon}(s)=s$ for all $s \in S$. 

\smallskip

(1) Suppose that $S$ is regular, and let $\alpha \in \Gamma$ and $s \in S_{\alpha}$ be arbitrary. Then there exists $t \in S$ such that $sts=s$, where, necessarily, $st \in E(S) \subseteq S_{\varepsilon}$. Therefore
\[(s* \tau_{\alpha}^{-1}(t))*s = (s\tau_{\alpha}\tau_{\alpha}^{-1}(t))*s = st\tau_{\varepsilon}(s) = sts = s,\]
and so $S^{\tau}$ is regular.

\smallskip

(2) Suppose that $S$ is an inverse semigroup. Then idempotents in $S$ commute--see, e.g.,~\cite[Theorem 5.1.1]{Howie}. Now, let $u,v \in E(S^{\tau})$. Then $u,v \in E(S)$, by Lemma~\ref{idempt-lem}, and hence
\[u*v = u \tau_{\varepsilon}(v) = uv = vu = v \tau_{\varepsilon}(u) = v*u.\]
Since $S^{\tau}$ is regular, by (1), and has commuting idempotents, it is an inverse semigroup, by~\cite[Theorem 5.1.1]{Howie}.
\end{proof}

Recall that the \emph{natural partial order} $\leq$ on an inverse semigroup $S$ is defined by $s\leq t$ ($s,t \in S$) if $s=tu$ (equivalently, $s=ut$) for some $u \in E(S)$. In particular, if $u,v \in E(S)$, then $u \leq v$ means that $u=uv=vu$. Also if $s \leq t$ for some $s,t \in S$, then $s^{-1} \leq t^{-1}$, and $rs \leq rt$, $sr \leq tr$ for all $r \in S$. See~\cite[\S 5.2]{Howie} for more details.

Recall also that two elements $s,t$ in an inverse semigroup $S$ are \emph{compatible} if $s^{-1}t, st^{-1} \in E(S)$, and \emph{orthogonal} if $s^{-1}t = 0 = st^{-1}$. It is easy to see that if $S$ is graded and $s,t \in S$ are compatible but not orthogonal, then $|s| = |t|$.

\begin{lemma} \label{npt}
Let $\, \Gamma$ be a group, $S$ a $\, \Gamma$-graded inverse semigroup, and $\tau = \{\tau_{\alpha} \mid \alpha \in \Gamma\}$ a Zhang twist on $S$, such that $\tau_{\varepsilon}(p)=p$ for all $p \in S$. Also let $\alpha \in \Gamma$, $s \in S_{\alpha}$, and $r,t \in S$.
\begin{enumerate}
\item[$(1)$] The inner inverse of $s$ in $S^{\tau}$ is $\tau_{\alpha}^{-1}(s^{-1})$.
\item[$(2)$] If $\tau$ is algebraic, then $\tau_{\alpha}^{-1}(s\inv)=\tau_{\alpha^{-1}}(s)\inv$.
\item[$(3)$] Letting $\leq$, respectively $\leq^{\tau}$, denote the natural partial order on $S$, respectively $S^{\tau}$, we have $\leq \, = \, \leq^{\tau}$.
\item[$(4)$] We have $r \leq t$ if and only if $\tau_{\alpha}(r) \leq \tau_{\alpha}(t)$.
\item[$(5)$] We have $\tau_{\alpha}(t)\tau_{\alpha}(t)^{-1} = \tau_{\alpha}(tt^{-1})$.
\item[$(6)$] Suppose that $t \in S_{\alpha}$. Then $s$ and $t$ are compatible, respectively orthogonal, in $S$ if and only if they are compatible, respectively orthogonal, in $S^{\tau}$.
\end{enumerate}
\end{lemma}

\begin{proof}
(1) As in the proof of Proposition~\ref{twist-reg}(1),
\[(s* \tau_{\alpha}^{-1}(s^{-1}))*s = (s\tau_{\alpha}\tau_{\alpha}^{-1}(s^{-1}))*s = ss^{-1}\tau_{\varepsilon}(s) = ss^{-1}s = s.\]
It follows that $s_{\tau}^{*}:=\tau_{\alpha}^{-1}(s^{-1})*s*\tau_{\alpha}^{-1}(s^{-1})$ is the unique inner inverse of $s$ in $S^{\tau}$. (In any inverse semigroup, if $xyx = x$ for some elements $x$ and $y$, then also $x(yxy)x=x$ and $(yxy)x(yxy)=yxy$.) Now, since $s^{-1} \in S_{\alpha^{-1}}$,
\[s_{\tau}^{*} = \tau_{\alpha}^{-1}(s^{-1})*(s\tau_{\alpha}(\tau_{\alpha}^{-1}(s^{-1})) = \tau_{\alpha}^{-1}(s^{-1}) * ss^{-1} = \tau_{\alpha}^{-1}(s^{-1})\tau_{\alpha^{-1}}(ss^{-1})\]
Applying (2) in Lemma~\ref{id-lem} gives
\[s_{\tau}^{*} = \tau_{\alpha}^{-1}(s^{-1}\tau_{\varepsilon}(ss^{-1})) = \tau_{\alpha}^{-1}(s^{-1}ss^{-1}) = \tau_{\alpha}^{-1}(s^{-1}).\]

\smallskip

(2) Suppose that $\tau$ is algebraic. Then $\tau_{\alpha}\tau_{\alpha^{-1}}(s^{-1}) = \tau_{\varepsilon}(s^{-1}) = s^{-1}$, and therefore $\tau_{\alpha}^{-1}(s^{-1}) = \tau_{\alpha^{-1}}(s^{-1})$. Since $\tau_{\alpha^{-1}}$ is, by hypothesis, a semigroup automorphism, it follows that $\tau_{\alpha}^{-1}(s^{-1}) = \tau_{\alpha^{-1}}(s)^{-1}$.

\smallskip

(3) Suppose that $s \leq t$. Then $s = ut$ for some $u \in E(S)$. Hence
\[u*t = u\tau_{\varepsilon}(t) = ut = s.\]
Since $u \in E(S^{\tau})$, by Lemma~\ref{idempt-lem}, we conclude that $s \leq^{\tau} t$. Thus $\leq$ is contained in $\leq^{\tau}$. The opposite inclusion follows from Proposition~\ref{inverse-twist}.

\smallskip

(4) Suppose that $r \leq t$, and write $r = ut$ for some $u \in E(S)$. Then
\[\tau_{\alpha}(r) = \tau_{\alpha}(ut) = \tau_{\alpha}(u\tau_{\varepsilon}(t)) = \tau_{\alpha}(u)\tau_{\alpha}(t).\]
Hence $\tau_{\alpha}(r) \leq \tau_{\alpha}(t)$, by Lemma~\ref{idempt-lem}.

Conversely, suppose that $\tau_{\alpha}(r) \leq \tau_{\alpha}(t)$. Then, using Lemma~\ref{idempt-lem}, we can write $\tau_{\alpha}(r) = \tau_{\alpha}(u)\tau_{\alpha}(t)$ for some $u \in E(S)$. Hence
\[\tau_{\alpha}(r) = \tau_{\alpha}(u)\tau_{\alpha}(t) = \tau_{\alpha}(u\tau_{\varepsilon}(t)) = \tau_{\alpha}(ut).\]
Since $\tau_{\alpha}$ is a bijection, we conclude that $r \leq t$.

\smallskip

(5) We have
\[\tau_{\alpha}(t)\tau_{\alpha}(t)^{-1} = \tau_{\alpha}(tt^{-1}t)\tau_{\alpha}(t)^{-1} = \tau_{\alpha}(tt^{-1}\tau_{\varepsilon}(t))\tau_{\alpha}(t)^{-1} = \tau_{\alpha}(tt^{-1})\tau_{\alpha}(t)\tau_{\alpha}(t)^{-1}.\]
Since $\tau_{\alpha}(t)\tau_{\alpha}(t)^{-1} \in E(S)$, we see that
\begin{equation} \label{ineq-eq}
\tau_{\alpha}(t)\tau_{\alpha}(t)^{-1} \leq \tau_{\alpha}(tt^{-1}).
\end{equation}

To show the opposite inequality we note that, by Proposition~\ref{inverse-twist}(1), $\tau^{-1} = \{\tau_{\alpha}^{-1} \mid \alpha \in \Gamma\}$ is a Zhang twist on $S^{\tau}$. Hence, substituting $\tau_{\alpha}^{-1}$ for $\tau_{\alpha}$ in equation~(\ref{ineq-eq}) gives
\[\tau_{\alpha}^{-1}(t)\tau_{\alpha}^{-1}(t)^{-1} \leq^{\tau} \tau_{\alpha}^{-1}(tt^{-1}).\]
Further substituting $\tau_{\alpha}(t)$ for $t$ we obtain
\[tt^{-1} = \tau_{\alpha}^{-1}(\tau_{\alpha}(t))\tau_{\alpha}^{-1}(\tau_{\alpha}(t))^{-1} \leq^{\tau} \tau_{\alpha}^{-1}(\tau_{\alpha}(t)\tau_{\alpha}(t)^{-1}),\]
and, by (3), also
\[tt^{-1} \leq \tau_{\alpha}^{-1}(\tau_{\alpha}(t)\tau_{\alpha}(t)^{-1}).\]
Taking $\tau_{\alpha}$ of both sides, by (4), we conclude that $\tau_{\alpha}(tt^{-1}) \leq \tau_{\alpha}(t)\tau_{\alpha}(t)^{-1}$.

\smallskip

(6) Let $s_{\tau}^{*}$ and $t_{\tau}^{*}$ denote the inner inverses of $s$ and $t$ in $S^{\tau}$, respectively. Then, by (1),
\[s*t_{\tau}^{*} = s\tau_{\alpha}(\tau_{\alpha}^{-1}(t^{-1})) = st^{-1}.\] 
Similarly,
\[\tau_{\alpha}(s_{\tau}^{*}*t) = \tau_{\alpha}(\tau_{\alpha}^{-1}(s^{-1})\tau_{\alpha^{-1}}(t)) = \tau_{\alpha}(\tau_{\alpha}^{-1}(s^{-1})) \tau_{\alpha\alpha^{-1}}(t) = s^{-1}\tau_{\varepsilon}(t) = s^{-1}t.\]
Thus, if $s$ and $t$ are compatible in $S$, then they are compatible in $S^{\tau}$, by Lemma~\ref{idempt-lem} and the hypothesis that $\tau_{\alpha}$ is a bijection. Likewise, if $s$ and $t$ are orthogonal in $S$, then they are orthogonal in $S^{\tau}$, by Lemma~\ref{twist-id-elt}(1). The converses of these statements follow from Proposition~\ref{inverse-twist}.
\end{proof}

\subsection{Boolean inverse semigroups} \label{bool-sec}

Our next goal is to show that graded-Boolean inverse $\land$-semigroups are preserved by twisting. We shall revisit these semigroups in Section~\ref{duality-sect}, in connection with noncommutative Stone duality. We begin by recalling the relevant terminology.

A partially ordered set $(X, \leq)$ is a \textit{lattice}, if for all $x, y\in X$ there exists in $X$ a meet (i.e., greatest lower bound) $x \wedge y$ and a join (i.e., least upper bound) $x\vee y$. A lattice $X$ is \textit{distributive} if 
\[x\wedge(y\vee z) = (x\wedge y)\vee (x\wedge z) \text{ and } x\vee(y\wedge z) = (x\vee y) \wedge (x\vee z),\] 
for all $x, y, z\in X$. A \emph{Boolean algebra} is a distributive lattice $X$, with top and bottom elements (denoted $1$ and $0$, respectively), together with a unary \emph{complement} operator $\lnot$, satisfying $x \lor \lnot x = 1$ and $x \land \lnot x = 0$, for all $x \in X$.

Let $(X, \leq)$ be a partially ordered set with bottom element $0$, and let $Y \subseteq X$. We set 
\[Y^{\uparrow} := \{x \in X \mid \exists y \in Y \text{ such that } y \leq x\},\]
and
\[Y^{\downarrow} := \{x \in X \mid \exists y \in Y \text{ such that } y \geq x\}.\]
If $Y = Y^{\downarrow}$, then $Y$ is an \emph{order ideal}. An order ideal of the form $\{x\}^{\downarrow}$ ($x \in X$) is \emph{principal}. We say that $Y$ is \emph{downward directed} if it is nonempty, and for all $x,y \in Y$ there exists $z\in X$ such that $z \leq x$ and $z \leq y$. If $Y$ is downward directed, $Y = Y^{\uparrow}$, and $0 \notin Y$, then $Y$ is a \emph{filter}. A maximal filter is called an \emph{ultrafilter}.

It is well-known that for any inverse semigroup $S$, the idempotent set $E(S)$ forms a lower semilattice (i.e., a partially ordered set such that every pair of elements has a meet) relative to the natural partial order $\leq$, where $e \land f = ef$ for all $e,f \in E(S)$. If every pair of elements in $S$ has a meet relative to $\leq$, then $S$ is called an \emph{inverse $\land$-semigroup}.

\begin{definition} \label{bool-def}
Let $\Gamma$ be a group, and $S$ a $\Gamma$-graded inverse semigroup. Then $S$ is $\Gamma$-\emph{graded-Boolean}, or simply \emph{graded-Boolean}, if it satisfies the following conditions.
\begin{enumerate}
\item[$(1)$] For all $\alpha \in \Gamma$, every orthogonal pair $s,t \in S_{\alpha}$ has a join $s\vee t$ in $S$.

\smallskip 

\item[$(2)$] $E(S)$ is a distributive lattice, and every nonzero principal order ideal of $E(S)$ is a Boolean algebra.
\end{enumerate}
\end{definition}

The above definition generalizes (non-graded) \emph{Boolean} inverse semigroups, first introduced by Lawson in~\cite{Lawson3}. See~\cite{RZ2} for more details.

\begin{proposition} \label{bool-prop}
Let $\, \Gamma$ be a group, $S$ a $\, \Gamma$-graded inverse semigroup, and $\tau = \{\tau_{\alpha} \mid \alpha \in \Gamma\}$ a Zhang twist on $S$. 
\begin{enumerate}
\item[$(1)$] $S$ is an inverse $\land$-semigroup if and only if $S^{\tau}$ is.

\item[$(2)$] $S$ is a $\, \Gamma$-graded-Boolean inverse semigroup if and only if $S^{\tau}$ is.
\end{enumerate}
\end{proposition}

\begin{proof}
First, note that, by Lemma~\ref{twist-mult} and Proposition~\ref{twist-reg}(2), $S^{\tau}$ is a $\Gamma$-graded inverse semigroup. Also, by Proposition~\ref{id-fix}, we may assume that $\tau_{\varepsilon}(s)=s$ for all $s \in S$. 

\smallskip

(1) This follows from Lemma~\ref{npt}(3).

\smallskip

(2) By Lemma~\ref{npt}(6), for each $\alpha \in \Gamma$ a pair of elements of $S_{\alpha} = S_{\alpha}^{\tau}$ is orthogonal in $S$ if and only if it is orthogonal in $S^{\tau}$. Thus such pairs have joins in $S$ if and only if they have joins in $S^{\tau}$, by Lemma~\ref{npt}(3). Also, by Lemma~\ref{idempt-lem} and Lemma~\ref{npt}(3), the lower semilattices $E(S)$ and $E(S^{\tau})$ are equal. Therefore $S$ satisfies the conditions in Definition~\ref{bool-def} if and only if $S^{\tau}$ does.
\end{proof}

\subsection{Properties preserved by twisting}

In this subsection we collect in one place various general semigroup properties, discussed above, that are preserved by twisting, and add some new ones.

Recall that a semigroup $S$ is a \emph{band} if $S = E(S)$, and an inverse semigroup $S$ is \emph{$0$-$E$-unitary} if $s \in E(S)$, for all $s \in S$ satisfying $u \leq s$ for some $u \in E(S)\setminus \{0\}$.

\begin{theorem} \label{twist-invar}
Let $\, \Gamma$ be a group, $S$ a $\, \Gamma$-graded semigroup, and $\tau = \{\tau_{\alpha} \mid \alpha \in \Gamma\}$ a Zhang twist on $S$. 
\begin{enumerate}
\item[$(1)$] $S$ is regular if and only if $S^{\tau}$ is.
\item[$(2)$] $S$ is an inverse semigroup if and only if $S^{\tau}$ is.
\item[$(3)$] $S$ is an inverse $\land$-semigroup if and only if $S^{\tau}$ is.
\item[$(4)$] $S$ is a $\, \Gamma$-graded-Boolean inverse semigroup if and only if $S^{\tau}$ is.
\item[$(5)$] $S$ is a $\, 0$-$E$-unitary inverse semigroup if and only if $S^{\tau}$ is.
\item[$(6)$] $S$ is a band if and only if $S^{\tau}$ is.
\item[$(7)$] $S$ is a monoid if and only if $S^{\tau}$ is.
\item[$(8)$] $S \setminus \{0\}$ is a group if and only if $S^{\tau} \setminus \{0\}$ is.
\end{enumerate}
\end{theorem}

\begin{proof}
Statements (1)--(4) are proved in Propositions~\ref{twist-reg} and~\ref{bool-prop}. For each of the remaining statements, by Proposition~\ref{inverse-twist}, it suffices to prove only the forward implication. Moreover, each of the properties in question is preserved by isomorphism, and so, by Proposition~\ref{id-fix}, we may assume that $\tau_{\varepsilon}(s)=s$ for all $s \in S$. 

\smallskip

(5) Suppose that $S$ is a $0$-$E$-unitary inverse semigroup. Then $S^{\tau}$ is an inverse semigroup, by $(2)$, and is $0$-$E$-unitary, by Lemma~\ref{idempt-lem} and Lemma~\ref{npt}(3).

\smallskip

(6) If $S$ is a band, then so is $S^{\tau}$, by Lemma~\ref{idempt-lem}.

\smallskip

(7) Suppose that $S$ is a monoid, and let $\alpha \in \Gamma$ and $s \in S_{\alpha}$ be arbitrary. By Lemma~\ref{twist-id-elt}(3), $1 = \tau_{\alpha}^{-1}(1)\tau_{\varepsilon}(1) = \tau_{\alpha}^{-1}(1)$, and hence $\tau_{\alpha}(1)=1$. Therefore
\[1*s = 1\cdot \tau_{\varepsilon}(s) = 1\cdot s = s = s\cdot 1 = s\tau_{\alpha}(1) = s *1,\]
and so $S^{\tau}$ is a monoid.

\smallskip

(8) Suppose that $S\setminus \{0\}$ is a group (with identity $1$). Then, by (7) and its proof, $S^{\tau}$ is a monoid, with identity $1$. Also, by (2), $S^{\tau}$ is an inverse semigroup. Since $E(S) = \{0,1\}$, we have $E(S^{\tau}) = \{0,1\}$, by Lemma~\ref{idempt-lem}. Now, let $s \in S \setminus \{0\}$, and let $s_{\tau}^{*}$ denote the inner inverse of $s$ in $S^{\tau}$. Then $s*s_{\tau}^{*}$ and $s_{\tau}^{*}*s$ are nonzero idempotents in $S^{\tau}$, and hence $s*s_{\tau}^{*} = 1 = s_{\tau}^{*}*s$. Therefore $S^{\tau}\setminus \{0\}$ is a group.
\end{proof}

\subsection{Corner skew monoids}\label{cornerskewsemigroup}

The construction of corner skew monoids given below appears to be new, and is inspired by that of corner skew Laurent polynomial rings in~\cite{araskew}, itself motivated by Paschke’s construction of certain $C^*$-algebras~\cite{paschke}. After defining these monoids, we show that the property of being a $\Z$-graded corner skew monoid is invariant under twisting. In Section~\ref{corneralgebra} we then extend this result to corner skew Laurent polynomial rings. A key motivation for these constructions is that Leavitt path algebras, of certain graphs, admit realizations as corner skew Laurent polynomial rings.

\begin{deff}
Let $S$ be a monoid, $S\langle x, y \rangle$ the monoid freely generated by $S$, $x$, and $y$ (i.e., the free product of $S$ and $\langle x, y \rangle$),  and $f : S\rightarrow uSu$ a semigroup isomorphism, where $u \in E(S)$ and $f(1)=u$. Define
\begin{align}\label{skewcornerla}
S[x,y,f] : = S\langle x, y \rangle / \langle yx=1, xy=u, xs=f(s)x, sy=yf(s) \mid s\in S\rangle. 
\end{align}
We refer to $uSu$ as a \emph{corner monoid}, and to $S[x,y,f]$ as a \emph{corner skew monoid}.
\end{deff}

Let $S[x,y,f]$ be a corner skew monoid and $\Gamma$ a group. Then assigning degree $\varepsilon$ to each nonzero element of $S$, and letting $|x|=|y|^{-1} \in \Gamma$, induces a $\Gamma$-grading on $S[x,y,f]$. In particular, letting $|x| :=1$ and $|y| :=-1$, gives $S[x,y,f]$ a natural $\mathbb{Z}$-grading. Note also that each element of $S[x,y,f]$ can be written in the form $sf^i(1)x^i$, for some $s \in S$ and $i\in \mathbb{N}$, and in the form $y^if^i(1)s$ ($s \in S$, $i\in \mathbb{N}$). 

The significance of this construction is that any $\mathbb{Z}$-graded monoid $S$, which has an element $t\in S_{-1}$ with a right inverse, is a corner skew monoid, as the following lemma shows. 

\begin{lemma}\label{chathapy}
Let $S$ be a $\mathbb{Z}$-graded monoid, with $t\in S_{-1}$ and $s\in S_1$ such that $ts=1$. Then $S\cong S_0[x,y,f]$ as graded semigroups, where 
\begin{align*}
f : S_0 &\longrightarrow st \, S_0 \, st\\
 r &\longmapsto s\, r\, t.
\end{align*}
\end{lemma}

\begin{proof}
We note that the degree-$0$ component $S_0$ is a submonoid of $S$. Since $ts=1$, the element $st$ is an idempotent. Hence the map $f$ is well-defined, and is an isomorphism of semigroups. The assignment $\phi:S_0[x,y,f]\rightarrow S$, where $\phi(x)=s$, $\phi(y)=t$, and $\phi(r)=r$ for all $r\in S_0$, satisfies the corner skew monoid relations (\ref{skewcornerla}), and induces a well-defined graded semigroup homomorphism. It is easy to see that $\phi$ is bijective in each homogeneous component, and thus $\phi$ is an isomorphism. 
\end{proof}

Let $E$ be the directed graph consisting of one vertex and one loop, and let $\SS(E)$ be the resulting graph inverse semigroup, as in Example~\ref{graphinsemigrp}. Then $\SS(E)$ can be identified with the \emph{bicyclic monoid} $S:= \langle x, y \mid yx =1 \rangle$ (with zero), where the $\mathbb{Z}$-grading is induced by letting $|x|:=1$ and $|y|: =-1$. By the previous lemma, $S \cong S_0[x,y,f]$, where $S_0=\{x^iy^i\mid i\in \mathbb{N}\}$, and $f: S_0\rightarrow xy\, S_0\, xy$ is given by $r \mapsto xry$. More generally, letting $E$ be the directed graph consisting of one vertex and $n \geq 1$ loops, the graph inverse semigroup $\SS(E)$, also known as the \emph{polycyclic monoid} in this case, can be realized as a skew corner monoid.

\begin{lemma}
Let $S[x,y,f]$ be a $\mathbb{Z}$-graded corner skew monoid, and set $u:= f(1)$.  
\begin{enumerate}
\item[$(1)$] $S$ is an inverse semigroup if and only if $S[x,y,f]$ is. 
\item[$(2)$] $S[x,y,f]$ is strongly $\mathbb{Z}$-graded if and only if $SuS=S$. 
 \end{enumerate}
\end{lemma}

\begin{proof}
(1) Suppose that $S$ is an inverse semigroup, and let $sf^i(1)x^i \in S[x,y,f]$ be an arbitrary element. An easy inductive argument shows that $x^iy^i = f^i(1)$ for all $i \in \N$. Since $f^i(1)$ is an idempotent, and since idempotents in $S$ commute, we have
\[(sf^i(1)x^i)(y^i f^i(1)s^{-1})(sf^i(1)x^i) = sf^i(1)s^{-1}sf^i(1)x^i = sf^i(1)x^i.\]
Likewise
\[(y^i f^i(1)s^{-1})(sf^i(1)x^i)(y^i f^i(1)s^{-1}) = y^i f^i(1)s^{-1},\]
and so $y^i f^i(1)s^{-1}$ is an inner inverse of $sf^i(1)x^i$. It is then easy to see that the only element of the form $y^j f^j(1)r$ that could play this role is $y^i f^i(1)s^{-1}$, and so $(sf^i(1)x^i)^{-1} = y^i f^i(1)s^{-1}$. Thus $S[x,y,f]$ is an inverse semigroup.

The converse is immediate. 

\smallskip

(2) Suppose that $S[x,y,f]$ is strongly graded. Since 
\[1\in S[x,y,f]_{0} = S[x,y,f]_1S[x,y,f]_{-1},\] 
$1=(rux)(yus)=rus$ for some $r,s \in S$. Thus for any $q\in S$ we have $q=q\cdot 1=qrus \in SuS$, and so $S=SuS$. 

The converse is also straightforward. 
\end{proof}

The following proposition shows that any twisting of a corner skew monoid, resulting from a graded automorphism (see Remark~\ref{aut-eg}), is again a corner skew monoid.

\begin{prop}\label{zhangskewsem}
Let $S[x,y,f]$ be a $\mathbb{Z}$-graded corner skew monoid, $\psi$ a graded automorphism of $S[x,y,f]$, and $\tau := \{\psi^n \mid n \in \Z\}$ the corresponding algebraic Zhang twist on $S[x,y,f]$. Then $S[x,y,f]^{\tau} \cong S[x,y,gf]$, via a graded isomorphism, where $g:=\psi|_S$.
\end{prop}

\begin{proof}
Since $\psi : S[x,y,f] \rightarrow S[x,y,f]$ is a graded automorphism, it restricts to a homogeneous component bijection $\psi_i : S[x,y,f]_i \rightarrow S[x,y,f]_i$ for each $i \in \mathbb{Z}$. In particular, $g := \psi_0 \ (=\psi|_S)$ is a semigroup automorphism of $S$. Since $\psi$ is graded, necessarily
\[\psi(x)=sx \text{ and } \psi(y)=yt,\]
for some $s\in Su$ and $t\in uS$, where $f: S\rightarrow uSu$ and $f(1)=u$. The corner skew monoid relations (\ref{skewcornerla}) then imply that
\begin{align}\label{gfyfhfh}
& ts = u, \\\notag 
& st = g(u), \\ \notag
& s\,fg(r) = gf(r)\,s \text{ for all } r \in S, \\
& t\,gf(r) = fg(r)\,t  \text{ for all } r \in S. \notag
\end{align}
(For example, the third equality can be obtained from $\psi(xr) = \psi(f(r)x)$.)

For convenience, let $T : = S[x,y,f]^{\tau}$. Since for $y \in T_{-1}$ and $sx \in T_1$, we have
\[y*(sx)= y\,\psi^{-1}(sx) = yx = 1,\] 
by Lemma~\ref{chathapy}, $T \cong S[x,y,h]$ as graded semigroups, where
\begin{align}\label{djjhdd}
h : S &\longrightarrow (sx*y)\,S \, (sx*y),\\
r &\longmapsto sx * r * y. \notag
\end{align}
Note that 
\[sx*y= sx\psi(y)=sxyt=sut=st=g(u).\]
Using the relations (\ref{gfyfhfh}) and (\ref{skewcornerla}), we have, for all $r\in S$,
\begin{multline*}
h(r) = sx*r*y  = (sx)*(ry) = sx\,\psi(ry) = sx\,g(r)\,\psi(y) = sx\,g(r)\,yt = s \, fg(r) \, xyt \\ 
= s \, fg(r) \, ut = s \, fg(r) \,t  = gf(r) \, st = gf(r)g(u) = gf(r) = g(u)gf(r)g(u).
\end{multline*}
So (\ref{djjhdd}) takes the form 
\begin{align*}
h : S &\longrightarrow g(u) S g(u)\\
r &\longmapsto gf(r).
\end{align*}
Thus $S[x,y,f]^{\tau} \cong S[x,y,gf]$, via a graded isomorphism. 
\end{proof}

%\textcolor{red}{what is $S[x,y,fg]$?}

\subsection{Twisted groups} \label{gp-sect}

Next we relate Zhang twists on groups to certain group extensions. In this subsection, when applied to groups, we interpret notions from above, such as \emph{graded}, without reference to a zero element, in the obvious way. 

Given a group $N$, we denote by $Z(N)$ the center of $N$, by $\mathrm{Aut}(N)$ the group of automorphisms of $N$, by $\mathrm{Inn}(N)$ the subgroup of inner automorphisms of $N$, and by $\mathrm{Out}(N): = \mathrm{Aut}(N)/\mathrm{Inn}(N)$ the group of \emph{outer} automorphisms of $N$. We have an exact sequence of groups 
\[1\to Z(N)\to N\xrightarrow{\iota} \mathrm{Aut}(N)\to \mathrm{Out}(N)\to 1,\] where $\iota_n(m) := nmn\inv$ for all $n,m \in N$.  As any automorphism of $N$ preserves $Z(N)$, and inner automorphisms act trivially, it follows that $\mathrm{Out}(N)$ acts on $Z(N)$ via $[\varphi](n) := \varphi(n)$, where $[\varphi]$ denotes the class of $\varphi\in \mathrm{Aut}(N)$.  

Given an exact sequence of groups
\[1\to N\to G\to \Gamma\to 1,\]
we can think of it as a $\Gamma$-grading of $G$, with $N = G_{\varepsilon}$. If $s\colon \Gamma\to G$ is any section satisfying $s(\varepsilon)=\varepsilon$, then we obtain a function, called a (\emph{normalized}) $2$-\emph{cocycle}, $c\colon \Gamma\times \Gamma\to N$ defined by 
\[s(\alpha)s(\beta) = c(\alpha,\beta)s(\alpha\beta).\]
Moreover, we have a map $\varphi\colon \Gamma\to \mathrm{Aut}(N)$, given by $\varphi_{\gamma}(n) := s(\gamma)ns(\gamma)^{-1}$, and it is immediate that 
\begin{equation}\label{eq:cocycle}
\varphi_{\alpha}\varphi_\beta = \iota_{c(\alpha,\beta)}\varphi_{\alpha\beta}
\end{equation}
for all $\alpha, \beta \in \Gamma$. Thus $\gamma\mapsto [\varphi_{\gamma}]$ is a homomorphism $\Phi\colon \Gamma\to \mathrm{Out}(N)$, i.e., an outer action of $\Gamma$ on $N$. This homomorphism is independent of the choice of $s$. Note that $\varphi$ induces a homomorphism $\Gamma\to \mathrm{Aut}(Z(N))$.

One verifies the following identities:
\begin{gather}\label{ga:cocycle.normalized}
c(\gamma,\varepsilon)=\varepsilon=c(\varepsilon,\gamma),\\ \label{ga:cocycle.id}
c(\alpha,\beta)c(\alpha\beta,\gamma) = \varphi_{\alpha}(c(\beta,\gamma))c(\alpha,\beta\gamma).
\end{gather}
As a group, $G\cong N\times \Gamma$, with multiplication 
\[(n,\alpha)(m,\beta) := (n\varphi_{\alpha}(m)c(\alpha,\beta),\alpha\beta).\] 
One has that $(\varepsilon,\varepsilon)$ is the identity, and 
\[(n,\gamma)\inv := (\varphi_\gamma\inv(n\inv c(\gamma,\gamma\inv)\inv),\gamma\inv).\] See~\cite[Chapter IV]{Brown} for more details.

Given a second extension
\[1\to N\to G'\to \Gamma\to 1\] 
giving rise to the same homomorphism $\Phi\colon \Gamma\to \mathrm{Out}(N)$, we can always choose a section $t$ so that $t(\gamma)nt(\gamma)\inv = \varphi_\gamma(n)$. Let $f\colon \Gamma\times \Gamma\to N$ be the corresponding $2$-cocycle. Using equation~\eqref{eq:cocycle} one can check that $fc\inv(\alpha,\beta) := f(\alpha,\beta)c(\alpha,\beta)\inv$ is a $2$-cocycle on $Z(N)$, viewed as a $\mathbb{Z}\Gamma$-module. We shall view $G'$ as $N\times \Gamma$, with product 
\[(n,\alpha)(m,\beta) := (n\varphi_{\alpha}(m)f(\alpha,\beta),\alpha\beta).\]
Our goal is to show that $G'$ is a twisting of $G$, where each $\tau_\gamma$ is the identity on $N$, and every Zhang twist of this sort gives an extension with the same outer action on $N$. Note that any twisting of a group is also a group, by Theorem~\ref{twist-invar}(8).

\begin{thm} \label{group-ext}
If \[1\to N\to G\to \Gamma\to 1\] is a group extension with outer action $\Phi\colon \Gamma\to \mathrm{Out}(N)$, and $\tau = \{\tau_{\alpha} \mid \alpha \in \Gamma\}$ is a Zhang twist on $G$ such that $\tau_{\alpha}|_N=1_N$ for all $\alpha\in \Gamma$, then 
\[1\to N\to G^{\tau}\to \Gamma\to 1\] 
is an extension with outer action $\Phi$.  

Conversely, if 
\[1\to N\to G'\to\Gamma\to 1\] 
is a group extension with outer action $\Phi$, then $G'\cong G^\tau$ for a Zhang twist $\tau = \{\tau_{\alpha} \mid \alpha \in \Gamma\}$ on $G$, with $\tau_{\alpha}|_N=1_N$ for all $\alpha\in \Gamma$.
\end{thm}

\begin{proof}
Let $s\colon \Gamma\to G$ be a section with $s(\varepsilon)=\varepsilon$, and let $\varphi\colon \Gamma\to \mathrm{Aut}(N)$ be as above, with $c$ the corresponding $2$-cocycle.   

First let $\tau = \{\tau_{\alpha} \mid \alpha \in \Gamma\}$ be a Zhang twist on $G$ such that $\tau_\alpha|_N=1_N$ for all $\alpha\in \Gamma$. Trivially 
\[1\to N\to G^{\tau}\to \Gamma\to 1\] 
is an extension. We check that it has outer action $\Phi$.  Indeed, if we use the same section $s$, then 
\[s(\gamma)\ast n\ast \tau_{\gamma}\inv (s(\gamma)\inv)= s(\gamma)\tau_{\gamma}(n)s(\gamma)\inv = s(\gamma)ns(\gamma)\inv=\varphi_\gamma(n),\] 
by our assumption on $\tau$, along with Lemmas~\ref{twist-id-elt}(4) and~\ref{npt}(1).

For the converse, let $f$ be the $2$-cocycle corresponding to $G'$, with respect to a section $t$, with $t(\varepsilon)=\varepsilon$ and $t(\gamma)nt(\gamma)\inv = \varphi_\gamma(n)$.  Viewing $G$ as having underlying set $N\times \Gamma$, as above, define 
\[\tau_\alpha(n,\beta) := (n\varphi_\alpha\inv(fc\inv(\alpha,\beta)),\beta).\]
Let us verify that $\tau$ is a Zhang twist. We shall use that $fc\inv\colon \Gamma\times \Gamma\to Z(N)$ is a $2$-cocycle.

We compute 
\begin{align*}
 \tau_\alpha((m,\beta)\tau_{\beta}(n,\gamma))&=\tau_{\alpha}((m,\beta)(n\varphi_\beta\inv(fc\inv(\beta,\gamma)),\gamma) \\ &= \tau_\alpha(m\varphi_\beta(n)fc\inv(\beta,\gamma)c(\beta,\gamma),\beta\gamma)\\ &= (m\varphi_\beta(n)fc\inv(\beta,\gamma)c(\beta,\gamma)\varphi_\alpha\inv(fc\inv(\alpha,\beta\gamma)),\beta\gamma)\\ &= (m\varphi_\beta(n)fc\inv(\beta,\gamma)\varphi_\alpha\inv(fc\inv(\alpha,\beta\gamma))c(\beta,\gamma),\beta\gamma),
\end{align*}
as $fc\inv$ takes values in $Z(N)$. On the other hand, using that $\varphi$ is an action on $Z(N)$ and that $fc\inv$ takes on values in the center, we have 
\begin{align*}
\tau_\alpha(m,\beta)\tau_{\alpha\beta}(n,\gamma) &= (m\varphi_\alpha\inv(fc\inv(\alpha,\beta)),\beta)(n\varphi_{\alpha\beta}\inv(fc\inv(\alpha\beta,\gamma)),\gamma) \\ &=(m\varphi_\alpha\inv(fc\inv(\alpha,\beta)),\beta)(n\varphi_\beta\inv\varphi_{\alpha}\inv(fc\inv(\alpha\beta,\gamma)),\gamma) \\ &= (m\varphi_\alpha\inv(fc\inv(\alpha,\beta))\varphi_\beta(n)\varphi_\alpha\inv(fc\inv(\alpha\beta,\gamma))c(\beta,\gamma),\beta\gamma) \\ &= (m\varphi_\beta(n)\varphi_\alpha\inv(fc\inv(\alpha,\beta))\varphi_\alpha\inv(fc\inv(\alpha\beta,\gamma))c(\beta,\gamma),\beta\gamma).
\end{align*}
Therefore, $\tau_\alpha((m,\beta)\tau_{\beta}(n,\gamma)) =\tau_\alpha(m,\beta)\tau_{\alpha\beta}(n,\gamma)$, by applying equation~\eqref{ga:cocycle.id} to $fc\inv$.
Note that, by equation~\eqref{ga:cocycle.normalized},
\[\tau_\alpha(n,\varepsilon) = (n\varphi_\alpha\inv(fc\inv(\alpha,\varepsilon)),\varepsilon) = (n,\varepsilon).\] 

Finally, we check that $G^\tau\cong G'$, where we view $G'$ as $N\times \Gamma$ with multiplication 
\[(m,\alpha)(n,\beta) := (m\varphi_\alpha(n)f(\alpha,\beta),\alpha\beta).\]  
Indeed, in $G^\tau$ we have 
\begin{align*}
(m,\alpha)\ast (n,\beta) &= (m,\alpha)\tau_\alpha(n,\beta) \\ &= (m,\alpha)(n\varphi_\alpha\inv(fc\inv(\alpha,\beta)),\beta) \\ &= (m\varphi_\alpha(n)f(\alpha,\beta)c(\alpha,\beta)\inv c(\alpha,\beta),\alpha\beta)
\\ &= (m\varphi_\alpha(n)f(\alpha,\beta),\alpha\beta)
\end{align*}
as required.
\end{proof}

We conclude this subsection with an example of a (non-algebraic) Zhang twist on a group, having the properties described in Theorem~\ref{group-ext}.

\begin{example} \label{gp-eg}
Let $i \geq 2$, $G := \Z_{i} \times \Z_{2}$, and $\Gamma := \Z_{2}$, viewed as additive groups. Then $G$ is $\Gamma$-graded by $\phi : G \to \Gamma$, where $\phi(m,n) := n$. Setting $N: = \Z_{i}$ gives the exact sequence 
\[0\to N\to G\to \Gamma\to 0.\]
Now, let $\tau_{0} : G \to G$ be the identity function, and define $\tau_{1} : G \to G$ by $\tau_{1}(m,0):= (m,0)$, and $\tau_{1}(m,1):= (m+1,1)$ for all $m \in \Z_{i}$. (So, in the notation of Theorem~\ref{group-ext}, $\tau_{0}|_N=1_N$ and $\tau_{1}|_N=1_N$.) Clearly, each of $\tau_{0}$ and $\tau_{1}$ is a graded bijection. It is easy to verify that
\[\tau_{\alpha}((m_1,n_1) + \tau_{\beta}((m_2,n_2))) = \tau_{\alpha}((m_1,n_1)) + \tau_{\alpha \beta}((m_2,n_2))\]
for all $\alpha, \beta \in \Gamma$, $(m_1,n_1) \in G_{\beta}$, and $(m_2,n_2) \in G$, by separately considering the cases where $n_1 = 0$, $n_2 = 0$, and $n_1=1=n_2$. Hence $\tau: = \{\tau_{0}, \tau_{1}\}$ is a Zhang twist on $G$. It is also easy to check that $G^{\tau}$ is generated by $(0,1)$, and is therefore isomorphic to $\Z_{2i}$. If $i \geq 3$, then $\tau_{1}$ is not a homomorphism, and so $\tau$ is not algebraic.
\end{example}

\subsection{Acts} 

In this subsection we develop a theory of twisted $S$-sets, analogous to that of twisted modules from~\cite{Zhang}, beginning with the necessary terminology.

Let $S$ be a semigroup and $X$ a set. Then $X$ is a \emph{right $S$-set}, or a \emph{right $S$-act}, if there is an action $X \times S \rightarrow X$ of $S$ on $X$, such that $(xs)t=x(st)$ for all $s,t \in S$ and $x \in X$. We say that a right $S$-set $X$ is \emph{unitary} (or \emph{unital}) if $XS=X$ (i.e, for all $x \in X$ there exist $y \in X$ and $s \in S$ such that $ys=x$). Also, a right $S$-set $X$ is \emph{pointed} if there exists a ``zero" element $0_X \in X$, such that $x0 = 0_X$ for all $x \in X$ (in which case, necessarily $0_Xs = 0_X$ for all $s \in S$). \emph{Left unitary pointed $S$-sets} are defined analogously. From now on, unless mentioned otherwise, all $S$-sets are assumed to be pointed. 

For right (respectively, left) $S$-sets $X$ and $Y$, a function $\phi:X\rightarrow Y$ is an $S$-\emph{map} if $\phi(xs)=\phi(x)s$ (respectively, $\phi(sx) = s\phi(x)$) for all $s\in S$ and $x\in X$. Note that if $\phi:X\rightarrow Y$ is an $S$-map, then $\phi(0_X) = 0_Y$.

Let us now assume that $S$ is $\Gamma$-graded, and let $X$ be a right $S$-set. Given a function $\phi : X \setminus \{0_X\} \to \Gamma$ such that $\phi (xs) = \phi (x)|s|$, for all $s \in S$ and $x \in X$ satisfying $xs \neq 0_X$, we say that $X$ is \emph{$\Gamma$-graded} (\emph{by $\phi$}). In this situation, for each $\alpha \in \Gamma$, we set $X_\alpha := \phi^{-1}(\alpha)\cup \{0_X\}$, and refer to each $x\in X_\alpha \setminus \{0_X\}$ as having \emph{degree} $\alpha$, written $|x|:=\alpha$. Equivalently, $X$ being \emph{$\Gamma$-graded} can be defined as $X = \bigcup_{\alpha \in \Gamma} X_\alpha$ for some $X_\alpha \subseteq X$, such that $X_{\alpha} S_{\beta} \subseteq X_{\alpha\beta}$ for all $\alpha,\beta \in \Gamma$, and $X_{\alpha} \cap X_{\beta} = \{0_X\}$ for all distinct $\alpha,\beta \in \Gamma$. The term \emph{$\Gamma$-graded} is defined analogously for left $S$-sets.

\begin{lemma} \label{set-twist-mult}
Let $\, \Gamma$ be a group, $S$ a $\, \Gamma$-graded semigroup, $\tau = \{\tau_{\alpha} \mid \alpha \in \Gamma\}$ a Zhang twist on $S$, and $X$ a unitary pointed right $S$-set $\, \Gamma$-graded by $\phi : X \setminus \{0_X\} \to \Gamma$. Define an action of $S^\tau$ on $X$ via 
\[x * s := x \tau_{\alpha}(s),\] 
for all $\alpha \in \Gamma$, $x \in X_{\alpha}$, and $s \in S$. With respect to this action $X$ is a unitary pointed right $S^{\tau}$-set, which is $\, \Gamma$-graded by $\phi$.
\end{lemma}

\begin{proof}
Let $\alpha, \beta \in \Gamma$, $x \in X_{\alpha}$, $s \in S_{\beta}$, and $t \in S$. Then
\[(x*s)*t = (x \tau_{\alpha}(s))\tau_{\alpha\beta}(t) = x (\tau_{\alpha}(s)\tau_{\alpha\beta}(t)) = x (\tau_{\alpha}(s\tau_{\beta}(t)) = x* (s*t),\]
and so $X$ is a right $S^{\tau}$-set. Since $X$ is unitary as a right $S$-set, and each $\tau_{\alpha}$ is a bijection, $X$ is also unitary as a right $S^{\tau}$-set. Since $\tau_{\alpha}(0) = 0$ for each $\alpha \in \Gamma$, by Lemma~\ref{twist-id-elt}(1), we have $x*0 = x \tau_{\alpha}(0) = 0_X$ for all $\alpha \in \Gamma$ and $x \in X_{\alpha}$, and so $X$ is pointed as a right $S^{\tau}$-set, with the same zero element. 

Finally, suppose that $x*s \neq 0_X$. Then $x \tau_{\alpha}(s) \neq 0_X$, and so
\[\phi(x*s) = \phi(x \tau_{\alpha}(s)) = \phi (x)|\tau_{\alpha}(s)| = \phi (x)|s|,\]
since $\tau_{\alpha}$ is graded. Thus $X$ is $\Gamma$-graded by $\phi$, as a right $S^{\tau}$-set.
\end{proof}

\begin{definition} \label{set-twist-def}
Let $\Gamma$ be a group, $S$ a $\Gamma$-graded semigroup, $\tau$ a Zhang twist on $S$, and $X$ a $\Gamma$-graded right $S$-set. We denote the right $S^{\tau}$-set constructed from $X$ in Lemma~\ref{set-twist-mult} by $X^{\tau}$, and refer to it as the \emph{twisting of $X$ by $\tau$}.
\end{definition}

For a semigroup $S$, we denote by $\Gc{S}$ the category whose objects are $\Gamma$-graded unitary (pointed) right $S$-sets, and whose morphisms are \emph{graded} $S$-maps, that is, $S$-maps $\phi : X \to Y$ such that $\phi(X_\alpha) \subseteq Y_\alpha$ for all $\alpha \in \Gamma$. We refer the reader to~\cite{Bergman} for general category-theoretic background, and to~\cite{RZ1} for more details about $\Gc{S}$.

Given a $\Gamma$-graded right $S$-set $X$, for all $\alpha, \beta \in \Gamma$ let $X(\alpha)_\beta :=X_{\alpha\beta}$. Then 
\begin{equation*}
X(\alpha): = \bigcup_{\beta \in \Gamma}X(\alpha)_\beta
\end{equation*}
is a $\Gamma$-graded right $S$-set, called the \emph{$\alpha$-shift} of $X$, where $X(\alpha) = X$ as $S$-sets, but the grading is ``shifted" by $\alpha$. (For a $\Gamma$-graded left $S$-set $X$, the \emph{$\alpha$-shift} is defined by setting $X(\alpha)_\beta : = X_{\beta\alpha}$ for all $\beta \in \Gamma$.) This construction leads to a \emph{shift functor} $\mathcal{T}_{\alpha}$ for each $\alpha \in \Gamma$, defined by
\begin{align*}
\mathcal{T}_{\alpha}: \Gc{S} &\longrightarrow \Gc{S} \\
X &\longmapsto X(\alpha),\notag
\end{align*}
which takes each morphism to itself. It is easy to see that each $\mathcal{T}_{\alpha}$ is an isomorphism, and that $\mathcal{T}_{\alpha} \mathcal{T}_{\beta}= \mathcal{T}_{\beta\alpha}$ for all $\alpha, \beta \in \Gamma$. 

The following result is an analogue of~\cite[Theorem 3.1]{Zhang}.

\begin{theorem} \label{cat-iso}
Let $\, \Gamma$ be a group, $S$ a $\, \Gamma$-graded semigroup, and $\tau = \{\tau_{\alpha} \mid \alpha \in \Gamma\}$ a Zhang twist on $S$. Letting $\mathcal{F}(X) := X^{\tau}$, as in Definition~\ref{set-twist-def}, for each object $X$ in $\, \Gc{S}$, and letting $\mathcal{F}(\psi) = \psi$ for each morphism $\psi$ in $\, \Gc{S}$, defines an isomorphism of categories $\mathcal{F} : \Gc{S} \rightarrow \Gc{S^{\tau}}$. Moreover, $\mathcal{F}(S(\alpha)) \cong S^{\tau}(\alpha)$ for all $\alpha \in \Gamma$.
\end{theorem}

\begin{proof}
By Lemma~\ref{set-twist-mult}, to conclude that $\mathcal{F}$ is well-defined, it suffices to check that given a morphism $\psi : X \rightarrow Y$ in $\Gc{S}$, the function $\psi : X^{\tau} \rightarrow Y^{\tau}$ is a graded $S^{\tau}$-map. Thus let $\alpha \in \Gamma$, $x \in X_{\alpha}$, and $s \in S$. Then, using the hypothesis that $\psi$ is a graded $S$-map,
\[\psi(x*s) = \psi(x \tau_{\alpha}(s)) = \psi(x) \tau_{\alpha}(s) = \psi(x) * s,\]
and so $\psi$ is an $S^{\tau}$-map. Since the gradings on $S$ and $S^{\tau}$ coincide (see Lemma~\ref{set-twist-mult}), $\psi : X^{\tau} \rightarrow Y^{\tau}$ is also graded.

Next, let $\tau^{-1}:= \{\tau_{\alpha}^{-1} \mid \alpha \in \Gamma\}$. Then, by Proposition~\ref{inverse-twist}, $\tau^{-1}$ is a Zhang twist on $S^{\tau}$, and $ (S^{\tau})^{\tau^{-1}} = S^{\tau\tau^{-1}} = S$. Thus we can define a functor $\mathcal{G} : \Gc{S^{\tau}} \rightarrow \Gc{S}$ by letting $\mathcal{G}(X) = X^{\tau^{-1}}$ for each object $X$ in $\Gc{S^{\tau}}$, and letting $\mathcal{G}(\psi) = \psi$ for each morphism $\psi$ in $\Gc{S^{\tau}}$. This functor is well-defined, by the argument in the previous paragraph. Clearly $\mathcal{G} \circ \mathcal{F}(\psi) = \psi$ for every morphism $\psi$ in $\Gc{S}$, and $\mathcal{F} \circ \mathcal{G}(\psi) = \psi$ for every morphism $\psi$ in $\Gc{S^{\tau}}$. Hence, to conclude that $\mathcal{F}$ is an isomorphism, it suffices to check that $\mathcal{G} \circ \mathcal{F}(X) = X$ for every object $X$ in $\Gc{S}$ (and then $\mathcal{F} \circ \mathcal{G}(X) = X$ for every object $X$ in $\Gc{S^{\tau}}$, by symmetry). Thus let $X$ be an object in $\Gc{S}$, $\alpha \in \Gamma$, $x \in X_{\alpha}$, and $s \in S$. Denoting the action of $S^{\tau}$ on $X^{\tau}$ by $*$, and the action of $S$ on $(X^{\tau})^{\tau^{-1}}$ by $\star$, we have
\[x \star s = x* \tau_{\alpha}^{-1}(s) = x \tau_{\alpha}\tau_{\alpha}^{-1}(s) = xs.\]
Therefore $X = (X^{\tau})^{\tau^{-1}}$ as right $S$-sets, and so $\mathcal{F}$ is an isomorphism.

To prove the final claim, let $\alpha \in \Gamma$ be arbitrary. Since $\tau_{\alpha^{-1}}$ restricts to a bijection $S_{\beta} \to S_{\beta}$ for each $\beta \in \Gamma$, it induces a graded bijection $\tau_{\alpha^{-1}} : S^{\tau}(\alpha) \to S(\alpha)^{\tau}$. So it suffices to show that $\tau_{\alpha^{-1}}$ is an $S^{\tau}$-map. Thus let $\beta \in \Gamma$, $s \in S_{\beta}$, and $t \in S$, and let $\star$ denote the $S^{\tau}$-action on $S(\alpha)^{\tau}$. Since the $S^{\tau}$-action $*$ on $S^{\tau}(\alpha)$ is the same as on $S^{\tau}$, and since $s \in S_{\alpha\alpha^{-1}\beta} = S(\alpha)_{\alpha^{-1}\beta}$, we have
\[\tau_{\alpha^{-1}}(s*t) = \tau_{\alpha^{-1}}(s\tau_{\beta}(t)) = \tau_{\alpha^{-1}}(s)\tau_{\alpha^{-1}\beta}(t) = \tau_{\alpha^{-1}}(s) \star t.\]
So $\tau_{\alpha^{-1}} : S^{\tau}(\alpha) \to \mathcal{F}(S(\alpha))$ is an isomorphism in $\Gc{S^{\tau}}$.
\end{proof}

Theorem~\ref{cat-iso} shows that $S(\alpha)^\tau \cong S^{\tau}(\alpha)$ as right $S^{\tau}$-sets. In general we can establish a certain equality between shifts of right $S^{\tau}$-sets. 

\begin{remark}
Let $\Gamma$ be a group, $S$ a $\Gamma$-graded semigroup, $\tau = \{\tau_{\alpha} \mid \alpha \in \Gamma\}$ a Zhang twist on $S$, and $X$ an object in $\Gc{S}$. For each $\alpha \in \Gamma$ define an action of $S^\tau$ on $X$ via 
\[x \star s := x\tau_{\alpha\beta}(s),\] 
for all $x \in X_{\beta}$ and $s \in S$. Then it is easy to see that with this action $X$ becomes a ($\Gamma$-graded unitary) right $S^\tau$-set, which we denote by $X^{\tau(\alpha)}$. 

Clearly $X^{\tau(\varepsilon)}=X^\tau$, as $\Gamma$-graded right $S^{\tau}$-sets. One can further check that 
\[X(\alpha)^{\tau(\beta)} = X^{\tau(\beta \alpha^{-1})} (\alpha),\] 
for all $\beta \in \Gamma$, as $\Gamma$-graded right $S^{\tau}$-sets.
\end{remark}

Our next goal is to prove the converse of Theorem~\ref{cat-iso}, appropriately understood. For that we require a couple of preliminary results.

\begin{proposition} \label{delta-semi}
Let $\, \Gamma$ be a group, let $S$ be a $\, \Gamma$-graded semigroup, and for all objects $X,Y$ in $\, \Gc{S}$ let $\, \mathrm{Hom}_{S}(X,Y)$ denote the set of morphisms $X \rightarrow Y$ in $\, \Gc{S}$. Also let 
\[\Delta(S) := \bigsqcup_{\alpha \in \Gamma} \mathrm{Hom}_{S}(S(\alpha^{-1}),S)\] 
denote the \emph{$0$-disjoint union} of the $\, \mathrm{Hom}_{S}(S(\alpha^{-1}),S)$. $($That is, $\bigsqcup_{\alpha \in \Gamma} \mathrm{Hom}_{S}(S(\alpha^{-1}),S)$ is the disjoint union of the sets $\, \mathrm{Hom}_{S}(S(\alpha^{-1}),S) \setminus \{0\}$, together with a single zero function.$)$ 

Then $\Delta(S)$ is a $\, \Gamma$-graded semigroup, where $\Delta(S)_{\alpha} := \mathrm{Hom}_{S}(S(\alpha^{-1}),S)$ for each $\alpha \in \Gamma$. Moreover, if $S$ is a monoid, then $S \cong \Delta(S)$, as $\, \Gamma$-graded semigroups.
\end{proposition}

\begin{proof}
Let $\alpha, \beta \in \Gamma$, $\phi \in \mathrm{Hom}_{S}(S(\alpha^{-1}),S)$, and $\psi \in \mathrm{Hom}_{S}(S(\beta^{-1}),S)$. Then, for all $\gamma \in \Gamma$, 
\[\psi(S((\alpha\beta)^{-1})_{\gamma}) = \psi(S(\beta^{-1})_{\alpha^{-1}\gamma}) \subseteq S_{\alpha^{-1}\gamma} = S(\alpha^{-1})_{\gamma}.\]
Thus we can view $\phi\psi$ as an element of $\mathrm{Hom}_{S}(S((\alpha\beta)^{-1}),S)$, upon interpreting $\psi$ as $\mathcal{T}_{\alpha^{-1}}(\psi)$. The fact that $\mathcal{T}_{\gamma^{-1}}\mathcal{T}_{\delta^{-1}} = \mathcal{T}_{(\gamma\delta)^{-1}}$, for all $\gamma, \delta \in \Gamma$, implies that this operation is associative. Moreover,
\[\Delta(S)_{\alpha}\Delta(S)_{\beta} = \mathrm{Hom}_{S}(S(\alpha^{-1}),S)\mathrm{Hom}_{S}(S(\beta^{-1}),S) \subseteq \mathrm{Hom}_{S}(S((\alpha\beta)^{-1}),S) = \Delta(S)_{\alpha\beta},\]
and so $\Delta(S)$ is a $\Gamma$-graded semigroup.

Now suppose that $S$ is a monoid, and let $\alpha \in \Gamma$ and $\phi \in \mathrm{Hom}_{S}(S(\alpha^{-1}),S)$. Then $\phi(1) \in S_{\alpha}$, since $1 \in S(\alpha^{-1})_{\alpha}$. Moreover, for all $s \in S(\alpha^{-1})$ we have $\phi(s) = \phi(1\cdot s) = \phi(1)s$. So $\phi$ is uniquely determined by $\phi(1)$, and conversely each element of $S$ determines one in $\mathrm{Hom}_{S}(S(\alpha^{-1}),S)$. Thus letting $f(\phi) := \phi(1)$ defines a graded bijection $f: \Delta(S) \rightarrow S$. Also, for all $\alpha, \beta \in \Gamma$, $\phi \in \mathrm{Hom}_{S}(S(\alpha^{-1}),S)$, and $\psi \in \mathrm{Hom}_{S}(S(\beta^{-1}),S)$ we have
\[f(\phi)f(\psi) = \phi(1)\psi(1) = \phi(1 \cdot \psi(1)) = \phi(\psi(1)) = f(\phi\psi),\]
since $\phi$ is an $S$-map. Hence $S \cong \Delta(S)$, as $\Gamma$-graded semigroups.
\end{proof}

The next observation is an analogue of~\cite[Proposition 2.8]{Zhang}. 

\begin{proposition} \label{iso-twist}
Let $\, \Gamma$ be a group, and let $S$ and $T$ be $\, \Gamma$-graded semigroups. Then the following are equivalent.
\begin{enumerate}
\item[$(1)$] $T \cong S^{\tau}$ as graded semigroups, for some Zhang twist $\tau$ on $S$.
\item[$(2)$] For each $\alpha \in \Gamma$ there is a graded bijection $\phi_{\alpha} : T \to S$, such that $\phi_{\alpha}(t_1t_2) = \phi_{\alpha}(t_1)\phi_{\alpha\beta}(t_2)$ for all $\alpha, \beta \in \Gamma$, $t_1 \in T_{\beta}$, and $t_2 \in T$.
\end{enumerate}
\end{proposition}

\begin{proof}
$(1) \Rightarrow (2)$ Suppose that $\psi : T \rightarrow S^{\tau}$ is a graded isomorphism, where $\tau = \{\tau_{\alpha} \mid \alpha \in \Gamma\}$. For each $\alpha \in \Gamma$ let $\phi_{\alpha} : = \tau_{\alpha}\psi$. Then for all $\alpha, \beta \in \Gamma$, $t_1 \in T_{\beta}$, and $t_2 \in T$,
\begin{multline*}
\phi_{\alpha}(t_1t_2) = \tau_{\alpha}\psi(t_1t_2) = \tau_{\alpha}(\psi(t_1)*\psi(t_2)) = \tau_{\alpha}(\psi(t_1)\tau_{\beta}\psi(t_2))\\
= \tau_{\alpha}(\psi(t_1))\tau_{\alpha\beta}\psi(t_2)) = \phi_{\alpha}(t_1)\phi_{\alpha\beta}(t_2).
\end{multline*}

\smallskip

$(2) \Rightarrow (1)$ Suppose that $\{\phi_{\alpha} \mid \alpha \in \Gamma\}$ satisfies the conditions in (2). Let $\alpha, \beta \in \Gamma$, $s_1 \in S_{\beta}$, and $s_2 \in S$, and write $s_1 = \phi_{\alpha}(t_1)$ and $s_2 = \phi_{\alpha\beta}(t_2)$ for some $t_1 \in T_{\beta}$ and $t_2 \in T$. Then 
\[\phi_{\alpha}^{-1}(s_1s_2) = \phi_{\alpha}^{-1}(\phi_{\alpha}(t_1)\phi_{\alpha\beta}(t_2))  = \phi_{\alpha}^{-1}(\phi_{\alpha}(t_1t_2)) = t_1t_2 = \phi_{\alpha}^{-1}(s_1)\phi_{\alpha\beta}^{-1}(s_2).\]
Now, for each $\alpha \in \Gamma$ let $\tau_{\alpha} := \phi_{\alpha}\phi_{\varepsilon}^{-1}$. Then for all $\alpha, \beta \in \Gamma$, $s_1 \in S_{\beta}$, and $s_2 \in S$, we have
\[\phi_{\varepsilon}^{-1}(s_1(\phi_{\beta}\phi_{\varepsilon}^{-1}(s_2))) = \phi_{\varepsilon}^{-1}(s_1)\phi_{\beta}^{-1}(\phi_{\beta}\phi_{\varepsilon}^{-1}(s_2)) = \phi_{\varepsilon}^{-1}(s_1)\phi_{\varepsilon}^{-1}(s_2),\]
and so
\begin{multline*}
\tau_{\alpha}(s_1\tau_{\beta}(s_2)) = \phi_{\alpha}\phi_{\varepsilon}^{-1}(s_1\phi_{\beta}\phi_{\varepsilon}^{-1}(s_2)) = \phi_{\alpha}(\phi_{\varepsilon}^{-1}(s_1)\phi_{\varepsilon}^{-1}(s_2))\\
= \phi_{\alpha}\phi_{\varepsilon}^{-1}(s_1)\phi_{\alpha\beta}\phi_{\varepsilon}^{-1}(s_2) = \tau_{\alpha}(s_1)\tau_{\alpha\beta}(s_2).
\end{multline*}
Thus $\tau : = \{\tau_{\alpha} \mid \alpha \in \Gamma\}$ is a Zhang twist on $S$.

Finally, for all $\alpha, \beta \in \Gamma$, $t_1 \in T_{\beta}$, and $t_2 \in T$,
\[\phi_{\varepsilon}(t_1t_2) = \phi_{\varepsilon}(t_1)\phi_{\beta}(t_2) = \phi_{\varepsilon}(t_1)\phi_{\beta}\phi_{\varepsilon}^{-1}\phi_{\varepsilon}(t_2) = \phi_{\varepsilon}(t_1)\tau_{\beta}\phi_{\varepsilon}(t_2) = \phi_{\varepsilon}(t_1)*\phi_{\varepsilon}(t_2),\]
and so $\phi_{\varepsilon}: T \rightarrow S^{\tau}$ is a graded isomorphism.
\end{proof}

We are ready for our main result about $\Gc{S}$, which is an analogue of~\cite[Theorem 3.4]{Zhang} 

\begin{theorem}\label{thmprehd6}
Let $\, \Gamma$ be a group, and let $S$ and $T$ be $\, \Gamma$-graded monoids. Then $T \cong S^{\tau}$ for some Zhang twist $\tau$ on $S$ if and only if there is an equivalence of categories $\mathcal{F} : \Gc{S} \rightarrow \Gc{T}$, such that $\mathcal{F}(S(\alpha)) \cong T(\alpha)$ for all $\alpha \in \Gamma$.
\end{theorem}

\begin{proof}
The forward implication follows from Theorem~\ref{cat-iso}. 

For the converse, suppose that $\mathcal{F} : \Gc{S} \rightarrow \Gc{T}$ is an equivalence satisfying $\mathcal{F}(S(\alpha)) \cong T(\alpha)$ for all $\alpha \in \Gamma$. By Proposition~\ref{delta-semi}, it suffices to show that $\Delta(T) \cong \Delta(S)^{\tau}$, as graded semigroups, for some Zhang twist $\tau$ on $\Delta(S)$. 

For each $\alpha \in \Gamma$ let $\psi_{\alpha} : \mathcal{F}(S(\alpha)) \rightarrow T(\alpha)$ denote the graded isomorphism provided by hypothesis. Then $f_{\alpha\beta} (\theta) := \psi_{\alpha^{-1}}^{-1}\theta \psi_{\beta^{-1}\alpha^{-1}}$ defines a graded bijection
\[f_{\alpha\beta} : \mathrm{Hom}_{T}(T(\beta^{-1}\alpha^{-1}),T(\alpha^{-1})) \rightarrow \mathrm{Hom}_{T}(\mathcal{F}(S(\beta^{-1}\alpha^{-1})),\mathcal{F}(S(\alpha^{-1})))\]
for each $\alpha \in \Gamma$, where $\mathrm{Hom}_{T}$ refers to morphism sets in the category $\Gc{T}$. Thus letting $\phi_{\alpha}(\rho) := \mathcal{T}_{\alpha}\mathcal{F}^{-1}f_{\alpha\beta}\mathcal{T}_{\alpha^{-1}}(\rho)$, for all $\alpha, \beta \in \Gamma$ and $\rho \in \mathrm{Hom}_{T}(T(\beta^{-1}), T)$, defines a graded bijection $\phi_{\alpha} : \Delta(T) \rightarrow \Delta(S)$--see the diagram below.
\[\xymatrix{\mathrm{Hom}_{T}(T(\beta^{-1}), T)\ar[d]^{\mathcal{T}_{\alpha^{-1}}} & & \mathrm{Hom}_{S}(S(\beta^{-1}), S)\\
 \mathrm{Hom}_{T}(T(\beta^{-1}\alpha^{-1}),T(\alpha^{-1}))\ar[rr]^{\mathcal{F}^{-1}f_{\alpha\beta}\ \ } & & \mathrm{Hom}_{S}(S(\beta^{-1}\alpha^{-1}),S(\alpha^{-1}))\ar[u]^{\mathcal{T}_{\a}}}\]
By Proposition~\ref{iso-twist}, to conclude the proof it suffices to check that $\phi_{\alpha}(\rho_1\rho_2) = \phi_{\alpha}(\rho_1)\phi_{\alpha\beta}(\rho_2)$ for all $\alpha, \beta, \gamma \in \Gamma$, $\rho_1 \in \Delta(T)_{\beta} = \mathrm{Hom}_{T}(T(\beta^{-1}),T)$, and $\rho_2 \in \Delta(T)_{\gamma} = \mathrm{Hom}_{T}(T(\gamma^{-1}),T)$. Taking such $\alpha$, $\beta$, $\gamma$, $\rho_1$, and $\rho_2$, we first recall that in the semigroup $\Delta(T)$ the product $\rho_1\rho_2$ is understood to be $\rho_1\mathcal{T}_{\beta^{-1}}(\rho_2)$ (see the proof of Proposition~\ref{delta-semi} for more details). Thus
\begin{align*}
\phi_{\alpha}(\rho_1\rho_2) & = \mathcal{T}_{\alpha}\mathcal{F}^{-1}f_{\alpha(\beta\gamma)}\mathcal{T}_{\alpha^{-1}}(\rho_1\mathcal{T}_{\beta^{-1}}(\rho_2)) \\
& = \mathcal{T}_{\alpha}\mathcal{F}^{-1}f_{\alpha(\beta\gamma)}(\mathcal{T}_{\alpha^{-1}}(\rho_1)\mathcal{T}_{\beta^{-1}\alpha^{-1}}(\rho_2))\\
& = \mathcal{T}_{\alpha}\mathcal{F}^{-1}(\psi_{\alpha^{-1}}^{-1}\mathcal{T}_{\alpha^{-1}}(\rho_1)\mathcal{T}_{\beta^{-1}\alpha^{-1}}(\rho_2) \psi_{\gamma^{-1}\beta^{-1}\alpha^{-1}})\\
& = \mathcal{T}_{\alpha}\mathcal{F}^{-1}(\psi_{\alpha^{-1}}^{-1}\mathcal{T}_{\alpha^{-1}}(\rho_1)\psi_{\beta^{-1}\alpha^{-1}})\mathcal{T}_{\alpha}\mathcal{F}^{-1}(\psi_{\beta^{-1}\alpha^{-1}}^{-1} \mathcal{T}_{\beta^{-1}\alpha^{-1}}(\rho_2) \psi_{\gamma^{-1}\beta^{-1}\alpha^{-1}}) \\
& = (\mathcal{T}_{\alpha}\mathcal{F}^{-1}f_{\alpha\beta}\mathcal{T}_{\alpha^{-1}}(\rho_1))(\mathcal{T}_{\alpha} \mathcal{T}_{\alpha\beta}^{-1} \mathcal{T}_{\alpha\beta} \mathcal{F}^{-1}(f_{(\alpha\beta)\gamma} \mathcal{T}_{\beta^{-1}\alpha^{-1}}(\rho_2)))\\
& = \phi_{\alpha}(\rho_1)\mathcal{T}_{\alpha} \mathcal{T}_{\beta^{-1}\alpha^{-1}}\phi_{\alpha\beta}(\rho_2)\\
& = \phi_{\alpha}(\rho_1)\mathcal{T}_{\beta^{-1}}\phi_{\alpha\beta}(\rho_2) \\
& = \phi_{\alpha}(\rho_1)\phi_{\alpha\beta}(\rho_2), 
\end{align*}
as desired.
\end{proof}

\section{Groupoids}\label{groupoidsec}

In this section we extend the Zhang twist from graded semigroups to graded
ample groupoids, and investigate its interaction with graded noncommutative Stone duality~\cite{Steinberg1,HazratLi2022,RZ2}.

\subsection{Graded ample groupoids} 

Recall that a \emph{groupoid} $\mathcal{G}$ is a small category where every morphism is invertible. The set of objects of $\mathcal{G}$, also known as the \emph{unit space}, is denoted by $\mathcal{G}^{0}$, and we identify these objects with the corresponding identity morphisms. For each morphism $a \in \mathcal{G}$, the object $\domr(a):=a^{-1}a$ is the \emph{domain} of $a$, and $\ran(a):=aa^{-1}$ is its \emph{range}. A pair $a,b\in \mathcal{G}$ is \emph{composable} if $\domr(a) = \ran(b)$. For all $U,V\subseteq \mathcal{G}$, we define
\begin{equation}\label{pofgtryhf}
UV:=\{ab \mid a \in U, b \in V, \domr(a)=\ran(b)\},
\end{equation}
and
\begin{equation}\label{pofgtryhf2}
U^{-1}:=\{a^{-1} \mid a \in U\}.
\end{equation}

A \textit{topological groupoid} is a groupoid equipped with a topology making inversion and composition continuous. A topological groupoid $\mathcal{G}$ is \textit{\'etale} if $\domr: \mathcal{G} \to \mathcal{G}^{0}$ is a local homeomorphism, where $\mathcal{G}^{0}$ is endowed with the subspace topology. An open subset $U$ of an \'etale groupoid $\mathcal{G}$ is called a \emph{local bisection}, or \emph{slice}, if the restrictions $\domr|_{U}$ and $\ran|_{U}$ are injective. An \'etale groupoid $\mathcal{G}$ is \emph{ample} if the compact local bisections form a base for its topology, and $\mathcal{G}^{0}$ is Hausdorff in the induced topology.

Given a topological groupoid $\mathcal{G}$ and a (discrete) group $\Gamma$, we say that $\mathcal{G}$ is $\Gamma$-\emph{graded} if there is a continuous (or locally constant) functor $\phi:\mathcal{G} \rightarrow \Gamma$. Setting $\mathcal{G}_\alpha:=\phi^{-1}(\alpha)$ for each $\alpha \in \Gamma$, we have $\mathcal{G}= \bigcup_{\alpha \in \G} \mathcal{G}_\alpha$, where $\mathcal{G}_\alpha \mathcal{G}_\beta \subseteq \mathcal{G}_{\alpha\beta}$ for all $\alpha,\beta \in \Gamma$, and $\mathcal{G}_\alpha \cap \mathcal{G}_\beta = \emptyset$ for all distinct $\alpha,\beta \in \Gamma$. As with semigroups, we say that each $a \in \mathcal{G}_\alpha$ has \emph{degree} $\alpha$, and write $|a| := \alpha$. Note that $\{\mathcal G_\alpha\mid \alpha \in \Gamma\}$ is a clopen partition of $\mathcal G$.

\subsection{Twisted ample groupoids}

The following is an analogue of Definition~\ref{twist-def} for ample groupoids.

\begin{definition} \label{groupoid-twist-def}
Let $\Gamma$ be a group and $\mathcal{G}$ a $\Gamma$-graded ample groupoid. A \emph{Zhang twisting system}, or a \emph{Zhang twist, $\tau$ on $\mathcal{G}$}, is a collection $\{\tau_{\alpha} \mid \alpha \in \Gamma\}$, where each $\tau_{\alpha} : \mathcal{G} \to \mathcal{G}$ is a graded homeomorphism (i.e., $\phi(\mathcal{G}_{\alpha}) \subseteq \mathcal{G}_{\alpha}$ for all $\alpha \in \Gamma$), such that the following conditions are satisfied for all $\alpha, \beta \in \Gamma$:
\begin{itemize}
\item[(T1)] $\tau_{\alpha} (\mathcal G^0)\subseteq \mathcal G^0$;
\item[(T2)] $\tau_{\alpha}|_{\mathcal G^0}\tau_{\beta}|_{\mathcal G^0} = \tau_{\alpha\beta}|_{\mathcal G^0}$;
\item[(T3)] $\domr\circ \tau_{\alpha}=\tau_{\alpha}\circ \domr$;
\item[(T4)] $\ran\circ \tau_{\alpha}=\tau_{\alpha}\circ \ran$;
\item[(T5)] if $a\in \mathcal G_{\beta}$ and $b \in \mathcal{G}$ are such that $\domr(a)=\ran(\tau_{\beta}(b))$, then $\tau_{\alpha}(a\tau_{\beta}(b)) = \tau_{\alpha}(a)\tau_{\alpha\beta}(b)$.
\end{itemize}
\end{definition}
 
\begin{remark} \label{dom-rem}
Note that (T5) makes sense because if $\domr(a)=\ran(\tau_{\beta}(b))$, then 
\[\ran(\tau_{\alpha\beta}(b)) = \tau_{\alpha\beta}(\ran(b)) = \tau_{\alpha}(\tau_{\beta}(\ran(b)) = \tau_{\alpha}(\domr(a)) = \domr(\tau_{\alpha}(a)),\] 
using (T1)--(T4). Conversely, if $a,b \in \mathcal{G}$ are such that $\domr(\tau_{\alpha}(a)) = \ran(\tau_{\alpha\beta}(b))$, then 
\[\tau_{\alpha}(\ran(\tau_{\beta}(b))) = \tau_{\alpha}(\tau_{\beta}(\ran(b))) = \tau_{\alpha\beta}(\ran(b)) = \ran(\tau_{\alpha\beta}(b)) = \domr(\tau_{\alpha}(a)) = \tau_{\alpha}(\domr(a)),\]
and so $\domr(a) = \ran(\tau_{\beta}(b))$, since $\tau_{\alpha}$ is bijective.
\end{remark}

For example, if $\mathcal{G}_{0}$ is a $\Gamma$-graded inverse semigroup obtained by adjoining a zero element to a discrete group $\mathcal{G}$, and $\tau = \{\tau_{\alpha} \mid \alpha \in \Gamma\}$ is any (semigroup) Zhang twist on $\mathcal{G}_{0}$ such that $\tau_{\varepsilon}(1) = 1$, then $\tau$ satisfies (T1)--T(5), by Lemma~\ref{twist-id-elt}(1) and Lemma~\ref{idempt-lem}. So restricting $\tau$ to $\mathcal{G}$ gives a Zhang twist, in the sense of Definition~\ref{groupoid-twist-def}. For another example, see Section~\ref{graph-groupoid-sect}, where we discuss graph groupoids.

As with semigroups, we shall use the Zhang twist to construct new groupoids. We begin with some preliminary observations.

\begin{lemma} \label{p:tau.eps}
Let $\, \Gamma$ be a group, $\mathcal{G}$ a $\Gamma$-graded ample groupoid, and $\tau = \{\tau_{\alpha} \mid \alpha \in \Gamma\}$ a Zhang twist on $\mathcal{G}$.
\begin{enumerate}
\item[$(1)$] $\tau_\varepsilon=1_{\mathcal G}$, the identity functor on $\mathcal{G}$.
\item[$(2)$] Each $\tau_\alpha$ restricts to a functor on $\mathcal{G}_\varepsilon$.
\end{enumerate}
\end{lemma}

\begin{proof}
(1) First note that (T1)--(T2) imply that on $\mathcal G^0$, we have $\tau_\varepsilon\tau_\varepsilon=\tau_\varepsilon$, and hence $\tau_\varepsilon|_{\mathcal G^0} = 1_{\mathcal G^0}$.  Now let $a\in \mathcal G$.  Then 
\[\tau_\varepsilon^2(a) = \tau_\varepsilon(aa\inv \tau_\varepsilon(a))=\tau_\varepsilon(aa\inv)\tau_\varepsilon(a) = aa\inv \tau_{\varepsilon}(a)=\tau_\varepsilon(a).\]  
Thus $\tau_\varepsilon(a) = a$.

\smallskip

(2) We compute for composable $a,b\in \mathcal G_\varepsilon$ that 
\[\tau_\alpha(ab) = \tau_\alpha(a\tau_\varepsilon(b)) = \tau_\alpha(a)\tau_\alpha(b),\] by (1) and (T5).
\end{proof}

Given a graded ample groupoid $\mathcal{G}$ and a Zhang twist $\tau = \{\tau_{\alpha} \mid \alpha \in \Gamma\}$ on $\mathcal{G}$, we define the \emph{twisting of $\mathcal{G}$ by $\tau$}, denoted $\mathcal G^\tau$, as follows. As a topological space $\mathcal G^\tau$ is $\mathcal G$.  The domain and range maps are given on $a\in \mathcal G_\alpha$ by $\ran^\tau(a) := \ran(a)$ and $\domr^\tau(a) := \tau_{\alpha\inv}(\domr(a))$. If $a\in \mathcal G_\alpha$ and $b\in \mathcal G_\beta$ with $\tau_{\alpha\inv}(\domr(a)) = \ran(b)$, then 
\[a\ast b := a\tau_{\alpha}(b).\] 
This makes sense because by (T1)--(T4), we have 
\[\ran(\tau_\alpha(b)) =\tau_\alpha(\ran(b))= \tau_\alpha(\tau_{\alpha\inv}(\domr(a))) = \domr(a).\]   

\begin{thm} \label{gp-twist-thrm}
Let $\, \Gamma$ be a group, $\mathcal{G}$ a $\, \Gamma$-graded ample groupoid, and $\tau = \{\tau_{\alpha} \mid \alpha \in \Gamma\}$ a Zhang twist on $\mathcal{G}$. Then $\mathcal G^\tau$ is a $\, \Gamma$-graded ample groupoid, with the same grading as $\mathcal{G}$.
\end{thm}

\begin{proof}
Multiplication is continuous being the composition 
\[(a,b)\mapsto (a,|a|,b)\mapsto (a,\tau_{|a|}(b))\mapsto a\tau_{|a|}(b)\] 
of continuous maps. 
    
If $a\in \mathcal G_\alpha$ and $b\in \mathcal G_\beta$, then $a\ast b=a\tau_\alpha(b)\in \mathcal G_{\alpha\beta}$, and so 
\[\domr^\tau(a\ast b) = \tau_{(\alpha\beta)\inv}(\domr(a\tau_\alpha(b))) = \tau_{(\alpha\beta)\inv}\tau_\alpha(\domr(b)) = \tau_{\beta\inv}(\domr(b)) = \domr^{\tau}(b)\] 
using (T1)--(T3).  Also $\ran^\tau(a\ast b) = \ran(a\tau_\alpha(b))=\ran^\tau(a)$.  
    
We next check associativity. Suppose that $a\in \mathcal G_\alpha$, $b\in \mathcal G_\beta$, and $c\in \mathcal G_\gamma$, where $a,b$ and $b,c$ are composable in $\mathcal G^\tau$. In particular, $\tau_{\beta\inv}(\domr(b))=\ran(c)$, and so $\domr(b) = \ran(\tau_\beta(c))$ by (T1)--(T4).  Then we compute
\begin{align*}
(a\ast b)\ast c &= (a\tau_\alpha(b))\tau_{\alpha\beta}(c)\\ 
a\ast(b\ast c) &= a\tau_\alpha(b\tau_\beta(c))
\end{align*}
and associativity follows from (T5).

Next we check that $\mathcal G^\tau$ is a category, by making sure $\ran^\tau(a)\ast a=a$ and $a\ast \domr^\tau(a)=a$.  So, let $a\in \mathcal G_\alpha$.  Then 
\[\ran^\tau(a)\ast a = \ran(a)\tau_\varepsilon(a) = \ran(a)a=a,\] 
by Lemma~\ref{p:tau.eps}(1).  On the other hand, 
\[a\ast \domr^{\tau}(a)= a \tau_\alpha(\tau_{\alpha\inv}(\domr(a))) = a\domr(a)=a,\] 
using (T1)--(T2).

Note that $\mathcal G^\tau\to \Gamma$, $a\mapsto |a|$, is a grading on $\mathcal G^\tau$, as $|a\ast b| = |a||\tau_{|a|}(b)| = |a||b|$.

Finally, we claim that the inverse of $a\in \mathcal G_\alpha$ in $\mathcal G^\tau$ is given by $\tau_{\alpha}\inv(a\inv)=\tau_{\alpha\inv}(a)\inv$.  Indeed, 
\[\ran^{\tau}(\tau_{\alpha}\inv(a\inv)) = \ran(\tau_{\alpha}\inv(a\inv)) = \tau_{\alpha}\inv(\domr(a)) = \tau_{\alpha\inv}(\domr(a)) = \domr^{\tau}(a).\] 
We then compute $a\ast \tau_{\alpha}\inv(a\inv) = aa\inv = \ran^{\tau}(a)$. Thus $\tau_{\alpha}\inv(a\inv)$ is a right inverse to $a$.  On the other hand, $\domr^{\tau}(\tau_{\alpha\inv}(a)\inv) = \ran(a)$, and we have 
\[\tau_{\alpha\inv}(a)\inv\ast a = \tau_{\alpha\inv}(a)\inv\tau_{\alpha\inv}(a) = \domr(\tau_{\alpha\inv}(a)) =\domr^{\tau}(a).\]    
In particular, it follows the inversion is continuous, as it is the composition $a\mapsto (|a|\inv,a)\mapsto \tau_{|a|\inv}(a)\inv$ of continuous maps.
\end{proof}

We record here the formula for the inverse in $\mathcal G^{\tau}$ (cf.\ Lemma~\ref{npt}).

\begin{cor}\label{c:inverse.gpd.twist}
The inverse in $\mathcal G^\tau$ is given on $a\in \mathcal G_{\alpha}$ by $\tau_{\alpha}\inv(a\inv)=\tau_{\alpha\inv}(a)\inv$.
\end{cor}

\begin{remark} \label{groupoid-alg-twist}
As with semigroups, we can construct \emph{algebraic} Zhang twists on a groupoid, through homomorphisms to the graded automorphism group. (See Remark~\ref{alg-eg}.) 

Specifically, for a $\Gamma$-graded ample groupoid $\mathcal{G}$, let $\mathrm{Aut}_{\mathrm{gr}}(\mathcal{G})$ denote the group of graded homeomorphisms of $\mathcal{G}$, that are also functors. Now, suppose that $\tau : \Gamma \rightarrow \mathrm{Aut}_{\mathrm{gr}}(\mathcal{G})$ is a homomorphism, and let $\tau_{\alpha} := \tau(\alpha)$ for each $\alpha \in \Gamma$. Then for all $\alpha, \beta \in \Gamma$, $a \in \mathcal{G}_{\beta}$, and $b \in \mathcal{G}$ such that $\domr(a)=\ran(\tau_\beta(b))$, we have 
\[\tau_\alpha(a\tau_\beta(b)) = \tau_\alpha(a)\tau_\alpha(\tau_\beta(b)) = \tau_\alpha(a)\tau_{\alpha\beta}(b),\]
which shows that $\{\tau_{\alpha} \mid \alpha \in \Gamma\}$ satisfies (T5). It is clear that (T1)--(T4) are satisfied as well.

If $\Gamma = \Z$, then, as in Remark~\ref{aut-eg}, we can construct an algebraic Zhang twist $\{\psi^n \mid n \in \Z\}$ on $\mathcal{G}$, from any graded automorphism $\psi : \mathcal{G} \to \mathcal{G}$.
\end{remark}

\subsection{Duality} \label{duality-sect}

In~\cite{Lawson3} Lawson established a duality between the categories Boolean inverse $\land$-semigroups and Hausdorff ample topological groupoids, with appropriate morphisms, which was generalized to graded versions thereof in~\cite{RZ2}. Results of this sort are referred to as \emph{(noncommutative) Stone dualities}, since they generalize Stone’s representation theorem for Boolean algebras~\cite{Stone}. Or next goal is to examine how the Zhang twist behaves under the aforementioned graded duality. We begin by recalling the relevant terminology.

Given an inverse semigroup $S$, let 
\[\GP(S):=\{U \mid U \text{~is an ultrafilter on~} (S, \leq)\},\] 
where $\leq$ is the natural partial order. (See Section~\ref{bool-sec} for the definition of \emph{ultrafilter}.) For all $U,V \in \GP(S)$ define 
\[U^{-1} := \{a^{-1} \mid a \in U\},\] 
\[UV := \{ab \mid a \in U, b \in V\},\] 
and 
\[U\cdot V := (UV)^{\uparrow}.\] 
According to~\cite[Section 2.1]{Lawson3}, $\GP(S)$ is a groupoid with respect to the operation $\cdot$, where $\GP(S)^{0} = \{U^{-1}\cdot U \mid U \in \GP(S)\}$. According to \cite[Lemma 5.6]{RZ2}, if $S$ is a $\Gamma$-graded-Boolean inverse $\land$-semigroup, graded via $\phi : S\setminus \{0\} \to \Gamma$, then $\GP(S)$ is a Hausdorff ample groupoid, graded via $\psi: \GP(S) \to \Gamma$, defined by $\psi(U) := \phi(a)$ for all $U \in \GP(S)$, where $a \in U$ is arbitrary. 

Given an ample groupoid $\mathcal{G}$, let
\[\SG(\mathcal{G}) :=\{U\mid U \text{~is a compact local bisection of~} \mathcal{G}\}.\]
Then $\SG(\mathcal{G})$ is an inverse semigroup under the operations given in (\ref{pofgtryhf}) and (\ref{pofgtryhf2}) (see~\cite[Proposition 2.2.4]{paterson1}), with $\emptyset$ as the zero element. Supposing that $\mathcal{G}$ is $\Gamma$-graded, for some group $\Gamma$, we say that a local bisection $U \subseteq \mathcal{G}$ is \emph{homogeneous} if $U \subseteq \mathcal{G}_\alpha$ for some $\alpha \in \Gamma$, and set
\begin{equation*}
\SG^{\gr}(\mathcal{G}) := \{U\mid U \text{~is a homogeneous compact local bisection of~} \mathcal{G}\}.
\end{equation*}
Since $\mathcal{G}_{\alpha^{-1}} = \mathcal{G}_{\alpha}^{-1}$ for each $\alpha \in \Gamma$, it follows that $\SG^{\gr}(\mathcal{G})$ is an inverse subsemigroup of $\SG(\mathcal{G})$. Moreover, defining $\phi :\SG^{\gr}(\mathcal{G}) \backslash \{\emptyset\} \rightarrow \Gamma$ by $\phi(U) := \alpha$, whenever $U\subseteq \mathcal{G}_\alpha$, turns $\SG^{\gr}(\mathcal{G})$ into a $\Gamma$-graded inverse semigroup. According to \cite[Lemma 5.1]{RZ2}, if $\mathcal{G}$ is a Hausdorff, then $\SG^{\gr}(\mathcal{G})$ is a $\Gamma$-graded-Boolean inverse $\land$-semigroup. 

The above constructions lead to a duality. Specifically, according to \cite[Lemma 5.11]{RZ2} (see also \cite[Theorem 3.3]{Steinberg2}), $\mathcal{G} \cong \GP(\SG^{\gr}(\mathcal{G}))$, for any $\Gamma$-graded Hausdorff ample groupoid $\mathcal{G}$, and according to~\cite[Lemma 5.12]{RZ2}, $S \cong \SG^{\gr}(\GP(S))$ for any $\Gamma$-graded-Boolean inverse $\land$-semigroup $S$. These mappings extend to a contravariant dual equivalence between the categories of relevant groupoids and semigroups, with appropriate morphisms--see \cite[Theorem 5.13]{RZ2}.

The rest of the subsection is devoted to showing how to pass Zhang twists between $S$ and $\GP(S)$, as well as between $\mathcal{G}$ and $\SG^{\gr}(\mathcal{G})$, in the above context. We begin with a simple observation.

\begin{lemma} \label{slice-lem}
Let $\, \Gamma$ be a group, $\mathcal{G}$ a $\, \Gamma$-graded ample groupoid, and $\tau = \{\tau_{\alpha} \mid \alpha \in \Gamma\}$ a Zhang twist on $\mathcal{G}$. Then $\, \SG^{\gr}(\mathcal{G}) = \SG^{\gr}(\mathcal{G}^{\tau})$ as sets.
\end{lemma}

\begin{proof}
Let $U\subseteq\mathcal{G}$. Since $\tau_{\alpha\inv}$ is a bijection, and $\domr^{\tau}|_U =\tau_{\alpha\inv}\circ \domr|_U$, we conclude that $\domr|_U$ is injective if and only if $\domr^{\tau}|_U$ is injective. Also, $\ran^\tau|_U = \ran|_U$, and hence either both are injective or neither. Since the topologies and gradings on $\mathcal{G}$ and $\mathcal{G}^{\tau}$ coincide (see Theorem~\ref{gp-twist-thrm}), we conclude that $\SG^{\gr}(\mathcal{G}) = \SG^{\gr}(\mathcal{G}^{\tau})$.
\end{proof}

The following notion will be convenient for relating Zhang twists on groupoids to those on corresponding semigroups, in the aforementioned duality.

\begin{definition}
Let $\Gamma$ be a group, $S$ a $\Gamma$-graded inverse semigroup, and $\tau = \{\tau_{\alpha} \mid \alpha \in \Gamma\}$ a Zhang twist on $S$. Then $\tau$ is \emph{semi-algebraic} if it satisfies the following conditions.
\begin{enumerate}
\item[$(1)$] For all $u \in E(S)$, we have $\tau_{\varepsilon}(u) = u$.
\item[$(2)$] For all $\alpha \in \Gamma$ and $s \in S$, we have $\tau_{\alpha}\inv(s\inv)= \tau_{\alpha\inv}(s)\inv$.
\item[$(3)$] The set $\{\tau_{\alpha}|_{E(S)} \mid \alpha \in \Gamma\}$ is an algebraic Zhang twist on $E(S)$.
\end{enumerate}
\end{definition}

By Lemma~\ref{idempt-lem}, any Zhang twist $\tau$ on $S$, satisfying (1) above, restricts to a Zhang twist on $E(S)$, which is a subsemigroup of $S$. If $\tau$ happens to be algebraic,
then it is semi-algebraic, since in that case $\tau_{\varepsilon}$ is necessarily the identity map, (2) is satisfied, by Lemma~\ref{npt}(2), and (3) clearly holds.

Let $\mathcal{G}$ be a $\Gamma$-graded group, and let $S:=\mathcal{G}_{0}$ be the $\Gamma$-graded inverse semigroup obtained by adjoining a zero element to $\mathcal{G}$. Then any Zhang twist $\tau = \{\tau_{\alpha} \mid \alpha \in \Gamma\}$ on $S$, such that $\tau_{\varepsilon}(1) = 1$, is semi-algebraic, by Lemma~\ref{twist-id-elt}(1) and Corollary~\ref{c:inverse.gpd.twist}, but not necessarily algebraic. A specific example of this sort can be constructed by adjoining a zero element to the group in Example~\ref{gp-eg}.

\begin{prop} \label{groupoid-semigp-prop}
Let $\, \Gamma$ be a group, $\mathcal{G}$ a $\, \Gamma$-graded ample groupoid, and $\tau = \{\tau_{\alpha} \mid \alpha \in \Gamma\}$ a Zhang twist on $\mathcal{G}$. For each $\alpha \in \Gamma$ define $\overline{\tau}_{\alpha} : \SG^{\gr}(\mathcal{G}) \to \SG^{\gr}(\mathcal{G})$ by $\overline{\tau}_{\alpha}(U) := \tau_{\alpha}(U)$. Then $\overline{\tau} :=\{\overline{\tau}_{\alpha} \mid \alpha \in \Gamma\}$ is a semi-algebraic Zhang twist on $\, \SG^{\gr}(\mathcal{G})$, and $\, \SG^{\gr}(\mathcal{G}^{\tau}) = \SG^{\gr}(\mathcal{G})^{\overline{\tau}}$ as semigroups.
\end{prop}

\begin{proof}
For each $\alpha \in \Gamma$, since $\tau_{\alpha}$ is a graded homeomorphism, it induces a bijection on the set of homogeneous compact open subsets of $\mathcal{G}$. It then follows from (T3)--(T4) that $\overline{\tau}_{\alpha} : \SG^{\gr}(\mathcal{G}) \to \SG^{\gr}(\mathcal{G})$ is a graded bijection. Now let $\alpha, \beta \in \Gamma$, $U \in \SG^{\gr}(\mathcal{G})_{\beta}$, and $V \in \SG^{\gr}(\mathcal{G})$. Then, using (T5) and Remark~\ref{dom-rem},
\begin{align*}
\overline{\tau}_{\alpha}(U\overline{\tau}_{\beta}(V)) & = \{\tau_{\alpha}(a\tau_{\beta}(b)) \mid a \in U, b \in V,  \domr(a)=\ran(\tau_{\beta}(b)) \big\}\\ 
& = \{\tau_{\alpha}(a)\tau_{\alpha\beta}(b) \mid a \in U, b \in V,  \domr(\tau_{\alpha}(a))=\ran(\tau_{\alpha\beta}(b)) \big\}\\ 
& = \overline{\tau}_{\alpha}(U)\overline{\tau}_{\alpha\beta}(V).
\end{align*}
Hence $\overline{\tau}$ is a Zhang twist on $\SG^{\gr}(\mathcal{G})$.

It is easy to see that $E(\SG^{\gr}(\mathcal{G}))$ consists of compact open subsets of $\mathcal{G}^{0}$. So, using (T1)--(T2), for any $U \in E(\SG^{\gr}(\mathcal{G}))$ we have $\tau_{\varepsilon} (\tau_{\varepsilon} (U)) = \tau_{\varepsilon}(U)$. Since $\tau_{\varepsilon}$ is a bijection, this means that $\overline{\tau}_{\varepsilon}(U) = \tau_{\varepsilon} (U) = U$. Hence $\overline{\tau}$ restricts to a Zhang twist on $E(\SG^{\gr}(\mathcal{G}))$, by Lemma~\ref{idempt-lem}. Also, using Corollary~\ref{c:inverse.gpd.twist},
\[\overline{\tau}_{\alpha\inv}(U)\inv = \{\tau_{\alpha\inv}(a)\inv \mid a \in U \big\}  = \{\tau_{\alpha}\inv(a) \mid a \in U\inv \big\} = \overline{\tau}_{\alpha}\inv(U\inv)\]
for all $\alpha \in \Gamma$ and $U \in \SG^{\gr}(\mathcal{G})$. Now, let $E: = E(\SG^{\gr}(\mathcal{G}))$, $\alpha \in \Gamma$, and $U, V \in E(\SG^{\gr}(\mathcal{G}))$. Then, using Lemma~\ref{p:tau.eps}(2), we have
\begin{align*}
\overline{\tau}_{\alpha}|_E(UV) & = \{\tau_{\alpha}(ab) \mid a \in U, b \in V,  \domr(a)=\ran(b) \big\}\\ 
& = \{\tau_{\alpha}(a)\tau_{\alpha}(b) \mid a \in U, b \in V,  \domr(\tau_{\alpha}(a))=\ran(\tau_{\alpha}(b)) \big\}\\ 
& = \overline{\tau}_{\alpha}|_E(U)\overline{\tau}_{\alpha}|_E(V),
\end{align*}
and so $\overline{\tau}_{\alpha}|_E$ is a semigroup homomorphism. Finally, using (T2), for all $\alpha, \beta \in \Gamma$, and $U \in E(\SG^{\gr}(\mathcal{G}))$, gives $\overline{\tau}_{\alpha}|_E(\overline{\tau}_{\beta}|_E(U)) = \overline{\tau}_{\alpha\beta}|_E(U)$. Hence $\{\tau_{\alpha}|_{E(S)} \mid \alpha \in \Gamma\}$ is an algebraic Zhang twist on $E(\SG^{\gr}(\mathcal{G}))$, and so $\overline{\tau}$ is a semi-algebraic Zhang twist on $\SG^{\gr}(\mathcal{G})$.

By Lemma~\ref{slice-lem}, $\SG^{\gr}(\mathcal{G}^{\tau}) = \SG^{\gr}(\mathcal{G}) = \SG^{\gr}(\mathcal{G})^{\overline{\tau}}$ as sets, and so for the final claim it suffices to show that the multiplication operations on $\SG^{\gr}(\mathcal{G}^{\tau})$ and $\SG^{\gr}(\mathcal{G})^{\overline{\tau}}$ coincide. Let $U,V \in \SG^{\gr}(\mathcal{G}^{\tau}) = \SG^{\gr}(\mathcal{G})^{\overline{\tau}}$, with $U \subseteq \mathcal{G}_{\alpha}$. Then in $\SG^{\gr}(\mathcal{G}^{\tau})$,
\begin{align*}
UV&= \{a*b \mid a \in U, b \in V, \domr^{\tau}(a) =\ran^{\tau}(b)\} \\ 
& = \{a\tau_{\alpha}(b) \mid a \in U, b \in V, \tau_{\alpha^{-1}}(\domr(a))=\ran(b)\}\\
& = \{a\tau_{\alpha}(b) \mid a \in U, b \in V, \domr(a)=\ran(\tau_{\alpha}(b))\},
\end{align*}
where the last equality follows from (T4) and the fact that $\tau_{\alpha^{-1}}|_{\mathcal G^0} = \tau_{\alpha}^{-1}|_{\mathcal G^0}$, by (T2). On the other hand, in $\SG^{\gr}(\mathcal{G})^{\overline{\tau}}$,
\[U*V = U\overline{\tau}_{\alpha}(V) = \{a\tau_{\alpha}(b) \mid a \in U, b \in V,  \domr(a)=\ran(\tau_{\alpha}(b)) \big\}.\]
Hence the multiplication operations on $\SG^{\gr}(\mathcal{G}^{\tau})$ and $\SG^{\gr}(\mathcal{G})^{\overline{\tau}}$ agree.
\end{proof}

\begin{proposition} \label{semigp-groupoid-prop}
Let $\, \Gamma$ be a group, $S$ a $\, \Gamma$-graded-Boolean inverse $\land$-semigroup, and $\tau = \{\tau_{\alpha} \mid \alpha \in \Gamma\}$ a semi-algebraic Zhang twist on $S$. For each $\alpha \in \Gamma$ define $\overline{\tau}_{\alpha} : \GP(S) \to \GP(S)$ by $\overline{\tau}_{\alpha}(U) := \tau_{\alpha}(U)$. Then $\overline{\tau} :=\{\overline{\tau}_{\alpha} \mid \alpha \in \Gamma\}$ is a Zhang twist on $\, \GP(S)$, and $\, \GP(S^{\tau}) = \GP(S)^{\overline{\tau}}$ as groupoids.
\end{proposition}

\begin{proof}
By Lemma~\ref{twist-id-elt}(4), $\tau_{\varepsilon}(s) = s$ for all $s \in S$, and so we can apply Lemma~\ref{npt}(4), to conclude that each $\tau_{\alpha}$ is an order-isomorphism on $S$, and therefore takes each ultrafilter on $S$ to an ultrafilter, as does $\tau_{\alpha}^{-1}$. Thus the function $\overline{\tau}_{\alpha} : \GP(S) \to \GP(S)$ is well-defined, for each $\alpha \in \Gamma$. Since each $\tau_{\alpha}$ is a graded bijection, so is $\overline{\tau}_{\alpha}$.

The topology on $\GP(S)$ is generated by basic open sets of the form 
\[\mathcal{Y}_s := \{U \in \GP(S) \mid s \in U\},\] 
where $s \in S$--see \cite[Lemma 5.6]{RZ2}. We have $\overline{\tau}_{\alpha}^{-1}(\mathcal{Y}_s) = \mathcal{Y}_{\tau_{\alpha}^{-1}(s)}$ and $\overline{\tau}_{\alpha}(\mathcal{Y}_s) = \mathcal{Y}_{\tau_{\alpha}(s)}$, for all $\alpha \in \Gamma$ and $s \in S$, and so $\overline{\tau}_{\alpha}$ is a homeomorphism.

By~\cite[Lemma 2.8(2)]{Lawson3}, an ultrafilter on $S$ belongs to $\GP(S)^0$ if and only if it contains an idempotent. Since each $\tau_{\alpha}$ preserves $E(S)$ (see Lemma~\ref{idempt-lem}), it follows that $\overline{\tau}_{\alpha}$ satisfies (T1).

To verify (T2), let $\alpha, \beta \in \Gamma$, $U \in \GP(S)$, and $r,s \in U$. Since $U$ is downward directed, we can find $t \in U$ such that $t\leq r$ and $t \leq s$. Write $t = ru$ for some $u \in E(S)$. Then $r^{-1}t = r^{-1}ru \in E(S)$, and so, given that $\tau$ is semi-algebraic, $\tau_{\alpha}(\tau_{\beta}(r^{-1}t)) = \tau_{\alpha\beta}(r^{-1}t)$. Since $r^{-1}t \leq r^{-1}s$, and $\tau_{\alpha\beta}$ is an order-isomorphism, we have $\tau_{\alpha}(\tau_{\beta}(r^{-1}t)) \leq \tau_{\alpha\beta}(r^{-1}s)$. Therefore 
\[\tau_{\alpha\beta}(r^{-1}s) \in \overline{\tau}_{\alpha}(\overline{\tau}_{\beta}(U^{-1}U))^{\uparrow} = \overline{\tau}_{\alpha}(\overline{\tau}_{\beta}(U^{-1} \cdot U)).\]
So $\overline{\tau}_{\alpha\beta}(U^{-1}U) \subseteq \overline{\tau}_{\alpha}(\overline{\tau}_{\beta}(U^{-1}\cdot U))$, which implies that $\overline{\tau}_{\alpha\beta}(U^{-1} \cdot U) \subseteq \overline{\tau}_{\alpha}(\overline{\tau}_{\beta}(U^{-1}\cdot U))$. Likewise $\tau_{\alpha\beta}(r^{-1}t) \leq \tau_{\alpha}(\tau_{\beta}(r^{-1}s))$ implies that $\tau_{\alpha}(\tau_{\beta}(r^{-1}s)) \in \overline{\tau}_{\alpha\beta}(U^{-1} \cdot U)$, and therefore $\overline{\tau}_{\alpha}(\overline{\tau}_{\beta}(U^{-1}\cdot U)) = \overline{\tau}_{\alpha\beta}(U^{-1} \cdot U)$.

For (T3), let $\alpha \in \Gamma$, $U \in \GP(S)$, and $r,s \in U$. Then, as above, there exists $t \in U$ such that $t\leq r$ and $t \leq s$. Hence also $t^{-1} \leq s^{-1}$, and therefore $t^{-1}t \leq t^{-1}r \leq s^{-1}r$. By hypothesis on $\tau$, Proposition~\ref{inverse-twist}(1), and Lemma~\ref{npt}(5),
\[\tau_{\alpha}(t)^{-1}\tau_{\alpha}(t) = \tau_{\alpha^{-1}}^{-1}(t^{-1})\tau_{\alpha^{-1}}^{-1}(t^{-1})^{-1} = \tau_{\alpha^{-1}}^{-1}(t^{-1}t) = \tau_{\alpha}(t^{-1}t)^{-1} = \tau_{\alpha}(t^{-1}t),\]
where the last equality follows from the fact that $\tau_{\alpha}(t^{-1}t) \in E(S)$, by Lemma~\ref{idempt-lem}. Thus
\[\tau_{\alpha}(s^{-1}r) \geq \tau_{\alpha}(t^{-1}t) = \tau_{\alpha}(t)^{-1}\tau_{\alpha}(t),\]
and so
\[\overline{\tau}_{\alpha} \circ \domr(U) = \overline{\tau}_{\alpha}(U^{-1} \cdot U) = \overline{\tau}_{\alpha}(U^{-1} U)^{\uparrow} \subseteq (\overline{\tau}_{\alpha}(U)^{-1} \overline{\tau}_{\alpha}(U))^{\uparrow} = \domr \circ \overline{\tau}_{\alpha} (U).\]
Likewise
\[\tau_{\alpha}(s)^{-1}\tau_{\alpha}(r) \geq \tau_{\alpha}(t)^{-1}\tau_{\alpha}(t) = \tau_{\alpha}(t^{-1}t),\] 
and hence
\[\domr \circ \overline{\tau}_{\alpha} (U) = (\overline{\tau}_{\alpha}(U)^{-1} \overline{\tau}_{\alpha}(U))^{\uparrow} \subseteq \overline{\tau}_{\alpha}(U^{-1} U)^{\uparrow} = \overline{\tau}_{\alpha} \circ \domr(U).\]
Thus $\domr \circ \overline{\tau}_{\alpha} = \overline{\tau}_{\alpha} \circ \domr$.

The verification of (T4) is entirely analogous to that of (T3), but using the identity $\tau_{\alpha}(tt^{-1}) = \tau_{\alpha}(t)\tau_{\alpha}(t)^{-1}$ (see Lemma~\ref{npt}(5)) in place of $\tau_{\alpha}(t^{-1}t) = \tau_{\alpha}(t)^{-1}\tau_{\alpha}(t)$.

To check (T5), let $\alpha, \beta \in \Gamma$, $U \in \GP(S)_{\beta}$, and $V \in \GP(S)$ be such that $\domr(U) = \ran(\tau_{\beta}(V))$. Then
\begin{align*}
\overline{\tau}_{\alpha}(U \cdot \overline{\tau}_{\beta}(V)) = \overline{\tau}_{\alpha}(U\overline{\tau}_{\beta}(V))^{\uparrow} & = \{\tau_{\alpha}(s \tau_{\beta}(t)) \mid s \in U, t \in V\}^{\uparrow}\\
& = \{\tau_{\alpha}(s) \tau_{\alpha\beta}(t) \mid s \in U, t \in V\}^{\uparrow}\\
& = (\overline{\tau}_{\alpha}(U)\overline{\tau}_{\alpha\beta}(V))^{\uparrow} = \overline{\tau}_{\alpha}(U) \cdot \overline{\tau}_{\alpha\beta}(V).
\end{align*}

Finally, by Lemma~\ref{npt}(3), the ultrafilters on $S$ and $S^{\tau}$ coincide, and so $ \GP(S^{\tau}) = \GP(S) = \GP(S)^{\overline{\tau}}$ as sets. Now, let $U,V \in \GP(S^{\tau}) = \GP(S)^{\overline{\tau}}$, and let $\alpha := |U|$ and $\beta := |V|$. Then in $\GP(S^{\tau})$, using Lemma~\ref{npt}(1) and the hypothesis on $\tau$,
\[\domr(U) = \{\tau_{\alpha^{-1}}(s)^{-1}*t \mid s,t \in U\}^{\uparrow} = \{\tau_{\alpha^{-1}}(s)^{-1}\tau_{\alpha^{-1}}(t) \mid s,t \in U\}^{\uparrow},\]
and
\[\ran(V) = \{p*\tau_{\beta}^{-1}(r^{-1}) \mid p,r \in V\}^{\uparrow} = \{p\tau_{\beta}(\tau_{\beta}^{-1}(r^{-1})) \mid p,r \in V\}^{\uparrow} = \{pr^{-1} \mid p,r \in V\}^{\uparrow},\]
while in $\GP(S)^{\overline{\tau}}$, using Corollary~\ref{c:inverse.gpd.twist},
\[\domr^{\overline{\tau}}(U) = \overline{\tau}_{\alpha^{-1}}(U)^{-1} * U = (\overline{\tau}_{\alpha^{-1}}(U)^{-1}\overline{\tau}_{\alpha^{-1}}(U))^{\uparrow} = \{\tau_{\alpha^{-1}}(s)^{-1}\tau_{\alpha^{-1}}(t) \mid s,t \in U\}^{\uparrow},\]
and
\[\ran^{\overline{\tau}}(V) = \{pr^{-1} \mid p,r \in V\}^{\uparrow}.\]
Thus the unit spaces of $\GP(S^{\tau})$ and $\GP(S)^{\overline{\tau}}$ coincide, and $U$ and $V$ are composable in $\GP(S^{\tau})$ if and only if they are composable in $\GP(S)^{\overline{\tau}}$. So to conclude the proof it suffices to suppose that $U$ and $V$ are composable, and show that their products in $\GP(S^{\tau})$ and $\GP(S)^{\overline{\tau}}$ agree. In $\GP(S^{\tau})$ we have
\[U \cdot V = \{s*t \mid s \in U, t \in V\}^{\uparrow} = \{s\tau_{\alpha}(t) \mid s \in U, t \in V\}^{\uparrow},\]
while in $\GP(S)^{\overline{\tau}}$,
\[U*V = U \cdot \overline{\tau}_{\alpha}(V) = \{s\tau_{\alpha}(t) \mid s \in U, t \in V\}^{\uparrow}.\]
Thus $U \cdot V = U*V$, as claimed.
\end{proof}

The last two propositions give us a complete characterization of how Zhang twists on graded Hausdorff ample groupoids behave under the graded noncommutative Stone duality described above.

\begin{theorem} \label{duality-thrm}
Let $\, \Gamma$ be a group, $S$ a $\, \Gamma$-graded-Boolean inverse $\land$-semigroup, and $\mathcal{G}$ a $\, \Gamma$-graded Hausdorff ample groupoid.
\begin{enumerate}
\item[$(1)$] The construction in Proposition~\ref{semigp-groupoid-prop} provides a one-to-one correspondence between semi-algebraic Zhang twists on $S$ and Zhang twists on $\, \GP(S)$.
\item[$(2)$] The construction in Proposition~\ref{groupoid-semigp-prop} provides a one-to-one correspondence between Zhang twists on $\mathcal{G}$ and semi-algebraic Zhang twists on $\, \SG^{\gr}(\mathcal{G})$.
\end{enumerate}
\end{theorem}

\begin{proof}
(1) According to \cite[Lemma 5.12]{RZ2}, $\nu_{S} : S \to \SG^{\gr}(\GP(S))$, defined by $\nu_{S}(s) := \mathcal{Y}_s$, where $\mathcal{Y}_s = \{U \in \GP(S) \mid s \in U\}$, is a graded isomorphism. So given a semi-algebraic Zhang twist $\tau = \{\tau_{\alpha} \mid \alpha \in \Gamma\}$ on $\SG^{\gr}(\GP(S))$, setting $\tau_{\alpha}^{\nu} : = \nu_{S}^{-1}\tau_{\alpha}\nu_{S}$ defines a graded bijection $\tau_{\alpha}^{\nu} : S \to S$ for each $\alpha \in \Gamma$. We claim that $\tau^{\nu}: = \{\tau_{\alpha}^{\nu} \mid \alpha \in \Gamma\}$ is a semi-algebraic Zhang twist on $S$.

For all $\alpha, \beta \in \Gamma$, $s \in S_{\beta}$, and $t \in S$, we have
\begin{multline*}
\tau_{\alpha}^{\nu}(s\tau_{\beta}^{\nu}(t)) = \nu_{S}^{-1}(\tau_{\alpha}\nu_{S}(s\nu_{S}^{-1}\tau_{\beta}\nu_{S}(t))) = \nu_{S}^{-1}(\tau_{\alpha}(\nu_{S}(s)\tau_{\beta}\nu_{S}(t))) \\
= \nu_{S}^{-1}(\tau_{\alpha}(\nu_{S}(s))\tau_{\alpha\beta}(\nu_{S}(t))) = \nu_{S}^{-1}(\tau_{\alpha}(\nu_{S}(s)))\nu_{S}^{-1}(\tau_{\alpha\beta}(\nu_{S}(t))) = \tau_{\alpha}^{\nu}(s)\tau_{\alpha\beta}^{\nu}(t),
\end{multline*}
and so $\tau^{\nu}$ is a Zhang twist on $S$. Since $\tau$ is semi-algebraic, 
\[\tau_{\varepsilon}^{\nu}(u) = \nu_{S}^{-1}\tau_{\varepsilon}\nu_{S}(u) = \nu_{S}^{-1}\tau_{\varepsilon}(\mathcal{Y}_u) = \nu_{S}^{-1}(\mathcal{Y}_u) = u\]
for all $u \in E(S)$, and
\[(\tau_{\alpha}^{\nu})\inv(s\inv) = \nu_{S}^{-1}\tau_{\alpha}\inv\nu_{S}(s\inv) = \nu_{S}^{-1}\tau_{\alpha}\inv(\nu_{S}(s)\inv) = \nu_{S}^{-1}(\tau_{\alpha\inv}(\nu_{S}(s))\inv)
= \tau_{\alpha\inv}^{\nu}(s)\inv\]
for all $\alpha \in \Gamma$ and $s \in S$. Clearly, the restriction of each $\tau_{\alpha}^{\nu}$ to $E(S)$ is a semigroup automorphism. Also,
\[\tau_{\alpha}^{\nu}\tau_{\beta}^{\nu}(u) = \nu_{S}^{-1}\tau_{\alpha}\nu_{S}\nu_{S}^{-1}\tau_{\beta}\nu_{S}(u) = \nu_{S}^{-1}\tau_{\alpha\beta}\nu_{S}(u) = \tau_{\alpha\beta}^{\nu}(u)\]
for all $\alpha,\beta \in \Gamma$ and $u \in E(S)$. Therefore $\tau^{\nu}$ is a semi-algebraic Zhang twist on $S$.

For each Zhang twist $\tau$ on $\GP(S)$, letting $\overline{\tau}$ be the corresponding semi-algebraic Zhang twist on $\SG^{\gr}(\GP(S))$, as in Proposition~\ref{groupoid-semigp-prop}, we then obtain a semi-algebraic Zhang twist $\varphi(\tau) := \overline{\tau}^{\nu}$ on $S$. In the opposite direction, for each semi-algebraic Zhang twist $\tau$ on $S$, letting $\overline{\tau}$ be as in Proposition~\ref{semigp-groupoid-prop}, we have a Zhang twist $\psi (\tau): = \overline{\tau}$ on $\GP(S)$. To conclude the proof we shall check that the functions $\varphi$ and $\psi$ are mutually inverse.

For any Zhang twist $\tau = \{\tau_{\alpha} \mid \alpha \in \Gamma\}$ on $\GP(S)$, $\alpha \in \Gamma$, and $U \in \GP(S)$, we have
\[(\psi \circ\varphi(\tau))_{\alpha}(U) = \nu_{S}^{-1}\overline{\tau}_{\alpha}\nu_{S}(U) = \nu_{S}^{-1}\overline{\tau}_{\alpha}(\{\mathcal{Y}_s \mid s \in U\}) = \nu_{S}^{-1}(\{\mathcal{Y}_{\tau_{\alpha}(s)} \mid s \in U\}) = \tau_{\alpha}(U),\]
and so $\psi \circ \varphi$ is the identity function. Likewise, for any semi-algebraic Zhang twist $\tau = \{\tau_{\alpha} \mid \alpha \in \Gamma\}$ on $S$, $\alpha \in \Gamma$, and $s \in S$,
\begin{multline*}
(\varphi \circ \psi(\tau))_{\alpha}(s) = \nu_{S}^{-1}\overline{\overline{\tau}}_{\alpha}\nu_{S}(s) = \nu_{S}^{-1}\overline{\overline{\tau}}_{\alpha}(\{U \in \GP(S) \mid s \in U\}) \\
= \nu_{S}^{-1}(\{\overline{\tau}_{\alpha}(U) \in \GP(S) \mid s \in U\}) = \nu_{S}^{-1}(\mathcal{Y}_{\tau_{\alpha}(s)}) = \tau_{\alpha}(s),
\end{multline*}
and so $\varphi \circ \psi $ is the identity function.

\smallskip

(2) According to \cite[Lemma 5.11]{RZ2}, $\mu_{\mathcal{G}} : \mathcal{G} \to \GP(\SG^{\gr}(\mathcal{G}))$, defined by $\mu_{\mathcal{G}}(a) := \mathcal{X}_a$, where $\mathcal{X}_a = \{U \in \SG^{\gr}(\mathcal{G}) \mid a \in U\}$, is a graded homeomorphism and functor. So given a Zhang twist $\tau = \{\tau_{\alpha} \mid \alpha \in \Gamma\}$ on $\GP(\SG^{\gr}(\mathcal{G}))$, setting $\tau_{\alpha}^{\mu} : = \mu_{\mathcal{G}}^{-1}\tau_{\alpha}\mu_{\mathcal{G}}$ defines a graded homeomorphism $\tau_{\alpha}^{\mu} : \mathcal{G} \to \mathcal{G}$ for each $\alpha \in \Gamma$. We claim that $\tau^{\mu}: = \{\tau_{\alpha}^{\mu} \mid \alpha \in \Gamma\}$ is a Zhang twist on $\mathcal{G}$.

It is clear that $\tau^{\mu}$ satisfies (T1) and (T2). Since $\mu_{\mathcal{G}}$ is a functor, it commutes with the operations $\domr$ and $\ran$, on the respective groupoids, and hence so does $\mu_{\mathcal{G}}^{-1}$. Since each $\tau_{\alpha}$ likewise commutes with $\domr$ and $\ran$ on $\GP(\SG^{\gr}(\mathcal{G}))$, it follows that $\tau^{\mu}$ satisfies (T3) and (T4). Checking that $\tau^{\mu}$ satisfies (T5) is entirely analogous to the corresponding verification for  $\tau^{\nu}$ in (1). Thus $\tau^{\mu}$ is a Zhang twist on $\mathcal{G}$.

For each semi-algebraic Zhang twist $\tau$ on $\SG^{\gr}(\mathcal{G})$, letting $\overline{\tau}$ be the corresponding semi-algebraic Zhang twist on $\GP(\SG^{\gr}(\mathcal{G}))$, as in Proposition~\ref{semigp-groupoid-prop}, we then obtain a Zhang twist $\varphi(\tau) := \overline{\tau}^{\mu}$ on $\mathcal{G}$. In the opposite direction, letting $\overline{\tau}$ be as in Proposition~\ref{groupoid-semigp-prop}, for each Zhang twist $\tau$ on $\mathcal{G}$, we have a semi-algebraic Zhang twist $\psi (\tau): = \overline{\tau}$ on $\SG^{\gr}(\mathcal{G})$. It remains to check that the functions $\varphi$ and $\psi$ are mutually inverse.

For any semi-algebraic Zhang twist $\tau = \{\tau_{\alpha} \mid \alpha \in \Gamma\}$ on $\SG^{\gr}(\mathcal{G})$, $\alpha \in \Gamma$, and $U \in \SG^{\gr}(\mathcal{G})$, we have
\[(\psi \circ \varphi(\tau))_{\alpha}(U) = \mu_{\mathcal{G}}^{-1}\overline{\tau}_{\alpha}\mu_{\mathcal{G}}(U) = \mu_{\mathcal{G}}^{-1}\overline{\tau}_{\alpha}(\{\mathcal{X}_a \mid a \in U\}) = \mu_{\mathcal{G}}^{-1}(\{\mathcal{X}_{\tau_{\alpha} (a)} \mid a \in U\}) = \tau_{\alpha}(U),\]
and so $\psi \circ \varphi$ is the identity function. Likewise, for any Zhang twist $\tau = \{\tau_{\alpha} \mid \alpha \in \Gamma\}$ on $\mathcal{G}$, $\alpha \in \Gamma$, and $a \in \mathcal{G}$,
\begin{multline*}
(\varphi \circ \psi(\tau))_{\alpha}(a) = \mu_{\mathcal{G}}^{-1}\overline{\overline{\tau}}_{\alpha}\mu_{\mathcal{G}}(a) = \mu_{\mathcal{G}}^{-1}\overline{\overline{\tau}}_{\alpha}(\{U \in \SG^{\gr}(\mathcal{G}) \mid a \in U\}) \\
= \mu_{\mathcal{G}}^{-1}(\{\overline{\tau}_{\alpha}(U) \in \SG^{\gr}(\mathcal{G}) \mid a \in U\}) = \mu_{\mathcal{G}}^{-1}(\mathcal{X}_{\tau_{\alpha}(a)}) = \tau_{\alpha}(a),
\end{multline*}
and so $\varphi \circ \psi $ is the identity function.
\end{proof}

\section{Algebras}\label{algebrasec}

We now study the algebraic consequences of semigroup and groupoid twisting, in the contexts of corner skew Laurent polynomial rings, Leavitt path algebras, and Steinberg algebras. We begin by recalling some basic terminology and facts about Zhang twists for algebras.

Let $K$ be a commutative unital ring, $\Gamma$ a group, $A$ a $\Gamma$-graded $K$-algebra, and $\{\tau_{\alpha} \mid \alpha \in \Gamma\}$ a Zhang twisting system $\tau$ on $A$ (see Example~\ref{ring-twist}). Define a binary operation $*$ on $A$ by letting 
\[a * b := a \tau_{\alpha}(b),\] 
for all $\alpha, \beta \in \Gamma$, $a \in A_{\alpha}$, and $b \in A_{\beta}$, and extending $K$-linearly to all of $A$. Then, according to \cite[Section 2]{Zhang}, $(A,*)$ is a $K$-algebra, with the same $\Gamma$-grading as $A$, and with $\tau_{\varepsilon}^{-1}(1)$ as the multiplicative identity. We denote the $K$-algebra $(A,*)$ by $A^{\tau}$, and refer to it as the \emph{(Zhang) twisting of $A$ by $\tau$}, or, more simply, as a \emph{twisted} algebra. (Note that in \cite{Zhang}, $\Gamma$ is allowed to be a monoid, rather than a group.) Moreover, the twisting system $\{\tau_{\alpha} \mid \alpha \in \Gamma\}$ is \emph{algebraic} if $\tau_{\alpha}$ is an algebra automorphism and $\tau_{\alpha}\tau_{\beta} = \tau_{\alpha\beta}$, for all $\alpha, \beta \in \Gamma$. As in Remark~\ref{aut-eg}, if $A$ is $\Z$-graded, and $\psi : A \rightarrow A$ is a graded automorphism, then $\{\psi^n \mid n \in \Z\}$ is an algebraic Zhang twist on $A$.

\subsection{Corner skew Laurent polynomial rings}\label{corneralgebra}

We extend the corner skew monoid construction, given in Section~\ref{cornerskewsemigroup}, to build corner skew Laurent polynomial rings, which are instances of the so-called \emph{fractional skew monoid rings} constructed in~\cite{araskew}. Although the treatment is largely similar to the monoid case, it will enable us to realize certain Leavitt path algebras, and twistings thereof, as corner skew Laurent polynomial rings.

\begin{deff}
Let $R$ be a unital ring, and $f : R\rightarrow uRu$ a ring isomorphism, where $u \in R$ is an idempotent and $f(1)=u$. Define
\begin{align}\label{skewcornerring}
R[x,y,f]: = R\langle x, y \rangle / \langle yx=1, xy=u, xr=f(r)x, ry=yf(r) \mid r\in R\rangle. 
\end{align}
The ring $uRu$ is called a \emph{corner ring}, and we refer to $R[x,y,f]$ as a \emph{corner skew Laurent polynomial ring}.
\end{deff}

Let $A: = R[x,y,f]$ be a corner skew Laurent polynomial ring and $\Gamma$ a group. Then letting $A_{\varepsilon}: = R$, and assigning the same degree to $x$ and $y^{-1}$, induces a $\Gamma$-grading on $A$. In particular, letting $|x| :=1$ and $|y| :=-1$, gives $A$ a natural $\mathbb{Z}$-grading. To describe this grading in more detail, let $u_0 :=1$, and $u_n :=f^n(u_0)$ for all $n \geq 1$. Then the elements of $A$ can be written in the form
\[y^jr_{-j} +y^{j-1}r_{-j+1}+\dots+yr_{-1}+r_0 +r_1x+\dots +r_ix^i,\]
where $r_{-n} \in u_n R$ and $r_n \in R u_n$, for all $n\geq 0$. Hence, as a $\mathbb{Z}$-graded ring, $A=\bigoplus_{i\in \mathbb{Z}}A_i$, where 
\begin{align*}
A_i& := Ru_ix^i \, \text{ for  } i>0,\\
A_i& :=y^iu_{-i}R \, \text{ for } i<0,\\
A_0& :=R. 
\end{align*} 
Clearly, if $u=1$ and $f$ is the identity map, then $R[x,y,f]$ reduces to the familiar Laurent polynomial ring  $R[x,x^{-1}]$. 

Similarly to Lemma~\ref{chathapy}, if in a $\mathbb{Z}$-graded ring $R$, we have $b\in R_{-1}$ and $a \in R_1$ such that $ba=1$, then $R \cong R_0[x,y,f]$, where $f:R_0\rightarrow ab \, R_0 \, ab$ is given by $r \mapsto arb$. Likewise, we have a ring-theoretic analogue of Proposition~\ref{zhangskewsem}.

\begin{prop}\label{zhangskewring}
Let $R[x,y,f]$ be a $\mathbb{Z}$-graded corner skew Laurent polynomial ring, $\psi$ a graded automorphism of $R[x,y,f]$, and $\tau := \{\psi^n \mid n \in \Z\}$ the corresponding algebraic Zhang twist on $R[x,y,f]$. Then $R[x,y,f]^{\tau} \cong R[x,y,gf]$, via a graded isomorphism, where $g:=\psi|_R$.
\end{prop}

\begin{proof}
The proof is similar to that of  Proposition~\ref{zhangskewsem}, and we omit it. 
\end{proof}

We next recall Leavitt path algebras, and interpret certain Zhang twistings of these rings as corner skew Laurent polynomial rings. Zhang twists on Leavitt path algebras have been previously studied in~\cite{NamSrivastavaVien2024, NamSrivastava2026}.

Let $K$ be a field (or, more generally, a commutative unital ring), $E$ a directed graph, and $\SS(E)$ the corresponding graph inverse semigroup (see Example~\ref{graphinsemigrp}). Then the contracted semigroup ring $C_K(E):= K[\SS(E)]$ (i.e., the ring resulting from identifying the zero of $\SS(E)$ with the zero of the semigroup ring $K\SS(E)$) is known as the \emph{Cohn path $K$-algebra of $E$}. Let $I$ be the ideal of $C_K(E)$ generated by all elements of the form $v-\sum_{e\in \so^{-1}(v)} ee^{-1},$ where $v \in E^0$ is a regular vertex (i.e., one that emits a finite but nonzero number of edges). Then $L_K(E):= C_K(E)/I$ is called the \emph{Leavitt path $K$-algebra of $E$}. The natural $\Z$-grading on $\SS(E)$ induces one on $C_K(E)$, which in turn induces a grading on $L_K(E)$. We refer the reader to \cite{Abrams-Monograph} for more details on Leavitt path algebras. 

Now let $K$ be a commutative unital ring, let $E$ a be finite graph with no sources (i.e., vertices that receive no edges), and write $E^0=\{v_1,\dots,v_n\}$. For each $1\leq i \leq n$,  choose $e_i \in E^1$ such that $\ra(e_i)=v_i$, and set $a :=e_1+\dots+e_n \in L_K(E)_1$. Then $b := e^{-1}_1+\dots+e^{-1}_n \in L_K(E)_{-1}$ is a left inverse of $a$. Hence $L_K(E) \cong L_K(E)_0[x,y,f]$, where 
\begin{align*}
f:L_K(E)_0 &\longrightarrow ab\,L_K(E)_0 \,ab\\
 r &\longmapsto a\, r\,b.
\end{align*}
Moreover, if $\psi$ is a graded automorphism of $L_K(E)$, and $\tau := \{\psi^n \mid n \in \Z\}$ is the corresponding algebraic Zhang twist on $L_K(E)$, then, by Proposition~\ref{zhangskewring}, we have 
\[L_K(E)^{\tau} \cong L_K(E)_0[x,y,f]^{\tau} \cong L_K(E)_0[x,y,gf],\]
where $g :=\psi |_{L_K(E)_0}$. 

\subsection{Steinberg algebras}

In this subsection we recall Steinberg algebras, and discuss how to extend a Zhang twist on the underlying groupoid of such an algebra to one on the ring.

If $X$ is a topological space with a base consisting of compact open sets, and $K$ is a commutative unital ring, then $KX$ denotes the $K$-submodule of $K^X$ spanned by the characteristic functions of compact Hausdorff open subsets of $X$. If $X$ is discrete, this is the free $K$-module on $X$; if $X$ is Hausdorff, this consists of the locally constant functions $f\colon X\to K$ with compact support. If $\mathcal G$ is an ample groupoid, then $K\mathcal G$ is a $K$-algebra (called the \emph{Steinberg $K$-algebra of $\mathcal G$}) with respect to the convolution product
\[fg(a) :=\sum_{bc=a}f(b)g(c).\]
For the construction of these algebras see \cite{Steinberg1}. 
If $\mathcal{G}$ is a $\Gamma$-graded ample groupoid, then $K\mathcal{G} =\bigoplus_{\gamma\in \Gamma}K\mathcal G_{\gamma}$ is a $\Gamma$-grading on $K\mathcal{G}$~\cite{HazratLi2022}. 

Let $K$ be a commutative unital ring, $\mathcal{G}$ a $\Gamma$-graded ample groupoid, and $\tau = \{\tau_{\alpha} \mid \alpha \in \Gamma\}$ a Zhang twist on $\mathcal{G}$. Then we can define a Zhang twist on $K\mathcal G$ by putting $\tau_{\alpha}(f) := f\circ \tau_{\alpha}\inv$ for $\alpha\in \Gamma$, abusing notation. We have a number of things to verify. If $U\subseteq \mathcal G_{\beta}$ is a compact Hausdorff open set, then $\tau_{\alpha}(U)\subseteq \mathcal G_{\beta}$ is also a compact Hausdorff open set, and $\tau_{\alpha}(1_U) = 1_{\tau_{\alpha}(U)}$, where $1_U$ denotes the characteristic function of $U$, and likewise for $1_{\tau_{\alpha}(U)}$. Thus $\tau_{\alpha}\colon K\mathcal{G} \to K\mathcal{G}$ is a $K$-linear graded bijection. We proceed to the difficult part of the verification that $\tau = \{\tau_{\alpha} \mid \alpha \in \Gamma\}$ is a Zhang twist on $K\mathcal{G}$.

\begin{lemma}\label{p:gpd.alg.twist}
Let $K$ be a commutative unital ring, $\, \Gamma$ a group, $\mathcal{G}$ a $\, \Gamma$-graded ample groupoid, and $\tau = \{\tau_{\alpha} \mid \alpha \in \Gamma\}$ a Zhang twist on $\mathcal{G}$. If $f\in K\mathcal G_{\beta}$ and $g\in K\mathcal G_{\gamma}$, then $\tau_{\alpha}(f\tau_{\beta}(g)) = \tau_{\alpha}(f)\tau_{\alpha\beta}(g)$ for all $\alpha, \beta, \gamma \in \Gamma$.
\end{lemma}

\begin{proof}
We first compute 
\[\tau_{\alpha}(f\tau_{\beta}(g))(a) = \sum_{bc=\tau_{\alpha}\inv(a)}f(b)g(\tau_{\beta}\inv(c))=\sum_{\tau_{\alpha}(bc)=a}f(b)g(\tau_{\beta}\inv(c)).\]
Putting $d:=\tau_{\beta}\inv(c)$, and using that $f(b)\neq 0$ implies $b\in \mathcal G_{\beta}$, we can rewrite this as
\[\sum_{\tau_{\alpha}(b\tau_{\beta}(d))=a}f(b)g(d)=\sum_{\tau_{\alpha}(b)\tau_{\alpha\beta}(d)=a}f(b)g(d).\] 
Now letting $x:=\tau_{\alpha}(b)$ and $y:= \tau_{\alpha\beta}(d)$, we can rewrite the last expression as
\[\sum_{xy=a}f(\tau_{\alpha}\inv(x))g(\tau_{\alpha\beta}\inv(y)) = \sum_{xy=a}\tau_{\alpha}(f)(x)\tau_{\alpha\beta}(g)(y) =(\tau_{\alpha}(f) \tau_{\alpha\beta}(g))(a),\]
as required.
\end{proof}

It turns out that groupoid twisting commutes with the Steinberg algebra construction.

\begin{thm} \label{groupoid-algebra-thrm}
Let $K$ be a commutative unital ring, $\, \Gamma$ a group, $\mathcal{G}$ a $\, \Gamma$-graded ample groupoid, and $\tau = \{\tau_{\alpha} \mid \alpha \in \Gamma\}$ a Zhang twist on $\mathcal{G}$. Let $\tau$ also denote the associated Zhang twist on $K\mathcal G$, given by $\tau_{\alpha}(f) := f\circ \tau_{\alpha}\inv$. Then $K\mathcal{G}^{\tau} = (K\mathcal{G})^{\tau}$ as $K$-algebras.
\end{thm}

\begin{proof}
By construction, $K\mathcal{G}^{\tau} = K\mathcal{G} = (K\mathcal{G})^{\tau}$ as $K$-modules.

Now let $f\in K\mathcal G_\beta$ and $g\in K\mathcal G_\gamma$.  Then their product in $K\mathcal G^\tau$ is given by
\[fg(a) = \sum_{b\ast c=a}f(b)g(c)=\sum_{b\tau_{\beta}(c)=a}f(b)g(c),\]
since $f(b)\neq 0$ implies $b\in \mathcal G_\beta$. Substituting $d:=\tau_\beta(c)$ yields
\[fg(a) = \sum_{bd=a}f(b)g(\tau_\beta\inv(d)) = \sum_{bd=a}f(b)\tau_\beta(g)(d) = (f\tau_\beta(g))(a)=f\ast g(a).\]
Hence $K\mathcal{G}^{\tau} = (K\mathcal{G})^{\tau}$ as $K$-algebras.
\end{proof}

\subsection{Algebraic systems on graph groupoids} \label{graph-groupoid-sect}

We conclude with an algebraic Zhang twist construction on a class of graded ample groupoids, connected to Leavitt path algebras.

Let $E=(E^0,E^1,\ra,\so)$ be a finite directed graph (see Example~\ref{graphinsemigrp} for the notation) without sinks (i.e., vertices that emit no edges). The \emph{infinite path space} $E^\infty$ consists of all infinite paths $e_1e_2\cdots$ in $E$. This is a compact Hausdorff and totally disconnected space, for the topology with a base of \emph{cylinder sets} \[xE^\infty := \{xy \mid y \in E^\infty, \ra(x)=\so(y)\},\] 
where $x$ is a (finite) path in $E$, and $\so$ is extended to $E^\infty$ is the obvious way. The \emph{shift map} $T\colon E^\infty\to E^\infty$ is the local homeomorphism given by $T(e_0e_1e_2\cdots) = e_1e_2\cdots$. An \emph{automorphism} of the shift space $(E^\infty,T)$ is a self-conjugacy, that is, a homeomorphism $h\colon E^\infty\to E^\infty$ such that $h\circ T=T\circ h$.

The \emph{graph groupoid} $\mathcal{G}_E$ is the Deaconu-Renault groupoid of this dynamical system. It has space of objects $E^\infty$, and morphisms
\[\mathcal{G}_E := \{(x,m-n,y)\in E^\infty\times \mathbb{Z}\times E^\infty\mid T^m(x)=T^n(y)\}. \]
The domain, range, and inverse are given by 
\begin{align*}
\domr(x,k,y) & := y,\\
\ran(x,k,y) & :=x,\\
(x,k,y)\inv & := (y,-k,x),
\end{align*}
respectively, and the product is
\[(x,k,y)(y,\ell, z) := (x,k+\ell,z).\]
The groupoid $\mathcal G_E$ is $\mathbb{Z}$-graded via $\phi :\mathcal G \rightarrow \mathbb{Z}$, defined by $\phi(x,k,y) := k$. This groupoid is ample, for the topology which has as a base of compact local bisections all sets of the form 
\[Z(x,y) := \{(xz,|x|-|y|,yz)\mid z\in \ra(x)E^\infty\},\]
where $x,y$ are paths in $E$ with $\ran(x)=\ran(y)$. It is readily checked that a larger base of compact local bisections is given by all sets of the form
\[Z(U,(m,n),V) := \{(x,m-n,y)\mid x\in U, y\in V, T^m(x)=T^n(y)\},\] 
where $U,V \subseteq E^\infty$ are compact open, $(m,n)\in \mathbb{N}\times \mathbb{N}$, $T^m(U)=T^n(V)$, and $T|_U$, $T|_V$ are injective. For example, $Z(x,y) = Z(xE^\infty,(|x|,|y|), yE^\infty)$.

It is known that every eventual conjugacy of edge shifts, in the sense of Matui and Matsumoto~\cite{MatsumotoMatui2014,CarlsenRuizSimsTomforde}, induces a graded isomorphism of graph groupoids (and conversely). In particular, every automorphism of the shift space $(E^\infty,T)$ induces a graded automorphism of $\mathcal{G}_E$, which can then be used to construct an algebraic Zhang twist, as in Remark~\ref{groupoid-alg-twist}. For completeness, we give the details of this case.

\begin{prop}\label{popfh77}
Let $E$ be a finite directed graph without sinks, and $h$ an automorphism of $\, (E^\infty,T)$.  Then $\tau\colon \mathcal G_E\to \mathcal G_E$, given by $\tau(x,k,y) := (h(x),k,h(y))$, is a graded automorphism.
\end{prop}

\begin{proof}
The map $\tau$ is well-defined. Indeed,   if $k=m-n$ and $T^m(x)=T^n(y)$, then 
\[T^m(h(x)) = h(T^m(x)) = h(T^n(y)) = T^n(h(y)),\] 
and so $\tau(x,m-n,y) = (h(x),m-n,h(y))$ is in $\mathcal G_E$.  Clearly $\tau$ preserves the degree. It is immediate from the definitions that $\tau$ is a functor, as 
\[\tau((x,k,y)(y,\ell,z)) = (h(x), k+\ell,h(z)) =(h(x),k,h(y))(h(y), \ell, h(z)) = \tau(x,k,y)\tau(y,\ell,z).\]
    
Note that $\tau\inv$ is induced by the automorphism $h\inv$, as $\tau\inv(x,k,y) = (h\inv(x),k,h\inv(y))$. Continuity follows from the observation that 
\[\tau\inv(Z(U,(m,n),V)) = Z(h\inv(U),(m,n),h\inv(V)).\]
Note that 
\[T^m(h\inv(U)) = h\inv (T^m(U)) = h\inv (T^n(V)) = T^n(h\inv(V)).\] 
Also $T|_{h\inv (U)}$ is injective, as $(T\circ h\inv)|_U = h\inv\circ T|_U$, and similarly for $V$.   
\end{proof}

Any graph automorphism of a directed graph $E$ induces a self-conjugacy in the obvious way.  %It is known that if $E$ has no multiple edges, then all automorphisms \RB{of edge shifts?} are induced by graph automorphisms. \RB{Reference..} 
However, typically there are many more, and the automorphism group of the shift space is usually infinite~\cite{BoyleFranksKitchens}. In the next example we construct shift space automorphisms that do not generally result from graph automorphisms.

\begin{example} \label{shift-eg}
Let $E=(E^0,E^1,\ra,\so)$ be a directed graph, where $E^0 = \{v\}$, and $E^1$ has at least $3$ elements. Let $x$ be a path in $E$, that, viewed as a word in the elements of $E^1$, does not overlap itself--that is, no nonempty prefix of $x$ is also a suffix. Let $\sigma\colon E^1 \to E^1$ be any permutation fixing the first and last letter of $x$. Define $h_{\sigma,x}\colon E^\infty\to E^\infty$ so that $h_{\sigma,x}(z)$ coincides with $z \in E^\infty$, except that the letter immediately before each occurrence of $x$ in $z$ is replaced by $\sigma$ of that letter. Because occurrences of $x$ do not overlap, and because consecutive occurrences $xx$ are fixed, as $\sigma$ fixes the last letter of $x$, it is easy to see that $h_{\sigma,x}$ is an automorphism, with the same order as $\sigma$. 
\end{example}

The $\mathcal{G}_E$ construction can be extended to arbitrary directed graphs $E$, using a more complicated topology, and it can be shown that $K\mathcal{G}_E \cong L_K(E)$ as $K$-algebras, for any commutative unital ring $K$ \cite[Theorem 2.7]{rigby}. Taking $E$ and $h_{\sigma,x}$ as in Example~\ref{shift-eg}, one can check that $h_{\sigma,x}$ induces a graded $K$-algebra automorphism $\psi$ of the Leavitt path algebra $L_K(E)$, via
\begin{align*}
e & \longmapsto \sigma(e)xx^{-1}+\sum_{|y|=|x|,y\neq x}eyy^{-1},\\
e^{-1} & \longmapsto xx^{-1}\sigma(e)^{-1}+\sum_{|y|=|x|,y\neq x}yy^{-1}e^{-1},
\end{align*}
for $e\in E^1$. %For instance, if $E^1 = \{e,f,g\}$ and $\sigma\colon E^1 \to E^1$ is the permutation that transposes $f$ and $g$, then $h_{\sigma,e}$ replaces $f$ and $g$ whenever they precede an $e$. The corresponding graded $K$-algebra automorphism of  $L_K(E)$ is given by 
%\begin{align*}
%e &\longmapsto e\\
%f &\longmapsto gee^{-1}+f(ff^{-1}+gg^{-1})\\
%g &\longmapsto fee^{-1}+g(ff^{-1}+gg^{-1}).
%\end{align*}
Now, combining Propositions~\ref{zhangskewring}, Theorem~\ref{groupoid-algebra-thrm}, and Proposition~\ref{popfh77} gives
\[K\mathcal G_E^\tau = (K\mathcal G_E)^\tau\cong L_K(E)^\tau\cong L_K(E)_0[x, y, gf],\]
where $f:L_K(E)_0 \rightarrow xx^{-1} L_K(E)_0xx^{-1}$ is defined by $a \mapsto xax^{-1}$, $g :=\psi |_{L_K(E)_0}$, and $\tau := \{\psi^n \mid n \in \Z\}$ is the associated algebraic Zhang twist.

%I'd like to do the general case, but I keep making silly mistakes.    Let $h\colon E^\infty\to E^\infty$ be an automorphism.  By the Curtis-Hedlund-Lyndon theorem 
%there is an integer $n$ called the anticipation and a function $H\colon E^{n+1}\to E$ such that \[h(e_0e_1\cdots) = H(e_0\cdots e_n)H(e_1\cdots e_{n+1})H(e_2\cdots e_{n+2})\cdots. \]
%Let $h\inv$ be the inverse of $h$ with the anticipation $k$. % We can always assume that $h,h\inv$ have the same anticipation $n$ by increasing the anticipation of one of them. 
%We write $H'\colon E^{k+1}\to E$ for the  function associated to $h\inv$.
%My current guess is the corresponding $\ast$-automorphism of $L_E$ should be given by 
%\[e\longmapsto \sum_{|u|=n+k}H(eH'(u_0\cdots u_k)H'(u_1\cdots u_{k+1})\cdots H'(u_{n-1}\cdots u_{k+n-1}))uu^*. \]
%This seems to give the same answer for the previous examples. 

\subsection*{Acknowledgments.} 
R.\ Hazrat acknowledges the Australian Research Council Discovery Project DP230103184. Part of this work was done during his visit to Benjamin Steinberg in New York in January 2026.  B.\ Steinberg was supported by the NSF grant DMS-2452324, a Simons Foundation Collaboration Grant, award number 849561, the Australian Research Council Grant DP230103184, and Marsden Fund Grant MFP-VUW2411.

We thank Be'eri Greenfeld for the suggestion to explore the Zhang twist in the setting of semigroups.

\end{document}